\documentstyle[twoside,12pt]{article}
\newtheorem{prop}{Proposition}[section]
\newtheorem{opr}[prop]{Definition}
\newtheorem{theo}[prop]{Theorem}

\newtheorem{conj}[prop]{Conjecture}
\newtheorem{rem}[prop]{Remark}
\newtheorem{coro}[prop]{Corollary}
\newtheorem{exam}[prop]{Example}

\newtheorem{lem}[prop]{Lemma}

\newcommand{\proof}{\noindent {\bf Proof. }}
\title{On the Hodge Structure of Projective
Hypersurfaces in Toric Varieties}

\author{Victor V. Batyrev\thanks{Supported by DFG,
Forschungsschwerpunkt Komplexe Mannigfaltigkeiten.}\\
 Universit\"at-GH-Essen, Fachbereich  6,  Mathematik \\
 Universit\"atsstr. 3, Postfach 10 37 64, D-4300 Essen 1 \\
Federal Republic of Germany \\
e-mail: matf$\emptyset \emptyset$@vm.hrz.uni-essen.de
\and
David A. Cox\thanks{Research supported by NSF grant DMS-9301161.} \\
Department of Mathematics and Computer Science \\
Amherst College\\
Amherst, MA 01002\\
e-mail: dac@cs.amherst.edu}

\begin{document}

\date{June 22, 1993}

\maketitle

\thispagestyle{empty}

The purpose of this paper is to explain one  extension of the ideas
of the Griffiths-Dolgachev-Steenbrink method for describing the
Hodge theory of smooth (resp.~quasi-smooth) hypersurfaces in
complex projective spaces (resp.~in weighted projective spaces). The
main idea of this method is the representation of the
Hodge components $H^{d-1-p,p}(X)$ in the middle cohomology group
of projective hypersurfaces
\[ X = \{ z \in {\bf P}^{d} : f(z)=0 \} \]
in ${\bf P}^d = {\rm Proj}\,
{\bf C} \lbrack z_1, \ldots, z_{d+1} \rbrack$
using homogeneous components of the quotient of the
polynomial ring ${\bf C} \lbrack z_1, \ldots, z_{d+1} \rbrack$
by the ideal $J(f) = \langle \partial f / \partial z_1, \ldots,
\partial f / \partial z_{d+1} \rangle$.  Basic references are
\cite{dolgach,griff1,pet.steen,steen}.

In this paper, we consider hypersurfaces $X$ in compact
$d$-dimensional toric varieties ${\bf P}_{\Sigma}$ associated with
complete rational polyhedral fan $\Sigma$ of simplicial cones ${\bf
R}^d$. According to the theory of toric varieties
\cite{dem,mam,dan1,oda}, ${\bf P}_{\Sigma}$ is defined by glueing
together of affine toric varieties ${\bf A}_{\sigma} = {\rm Spec}\,
{\bf C} \lbrack \check{\sigma} \cap {\bf Z}^d
\rbrack$ $(\sigma \in \Sigma)$ where $\check{\sigma}$ denotes the dual
to $\sigma$ cone. Weighted projective spaces are examples of toric varieties.

M. Audin \cite{aud} first noticed that there exists another approach to
the definition of the toric variety ${\bf P}_{\Sigma}$.  This definition
bases on the
representation of ${\bf P}_{\Sigma}$
as a  quotient of some  Zariski open subset $U(\Sigma)$
in an affine space ${\bf A}^n$ by a linear diagonal action of
some algebraic subgroup ${\bf D}(\Sigma) \subset ({\bf C}^*)^n$.
The group of characters of ${\bf D}(\Sigma)$ is isomorphic to
the group of classes $Cl ({\Sigma})$ of divisors on ${\bf P}_{\Sigma}$
modulo the rational equivalence.
The dimenson $n$ of the open set $U(\Sigma)$ equals the number of
1-dimensional cones in the fan $\Sigma$, and the dimension of
${\bf D}(\Sigma)$ equals $n-d$, the rank of the Picard group of
${\bf P}_{\Sigma}$.  In particular, if ${\bf P}_{\Sigma}$ is smooth, then
$U(\Sigma)$ is the universal torsor over ${\bf P}_{\Sigma}$ (see
\cite{man.zfas}) and
${\bf D}(\Sigma)$ is the torus of Neron-Severi.

	The codimension of the complement
\[Z(\Sigma) = {\bf A}^n \setminus U(\Sigma)\]
is at least 2.  So the ring of regular algebraic functions on
$U(\Sigma)$ is isomorphic to the polynomial ring
\[S(\Sigma) = {\bf C}[z_1, \ldots, z_n].\]
The action ${\bf D}(\Sigma)$ on $U(\Sigma)$ induces a canonical
grading of the ring $S(\Sigma)$ by elements of $Cl(\Sigma)$, i.e., by
characters of ${\bf D}(\Sigma)$.  In the paper \cite{cox} of the
second author, the polynomial ring $S(\Sigma)$ together with the
$Cl(\Sigma)$-grading is called the {\em homogeneous coordinate ring of
the toric variety} ${\bf P}_{\Sigma}$.  One nice feature of this ring
is that a hypersurface $X \subset {\bf P}_\Sigma$ has a defining
equation $f = 0$ for some $f \in S(\Sigma)_\beta$.  Here,
$S(\Sigma)_\beta$ is is the graded piece of $S(\Sigma)$ in degree
$\beta$, and $\beta$ is the divisor class of $X$ in $Cl(\Sigma)$.

	Let us describe the contents of the paper in more detail:
\medskip

{\em Sections 1 and 2.}  We establish notation and review the
construction of a simplicial toric variety ${\bf P}_\Sigma$ as a
quotient $U(\Sigma)/{\bf D}(\Sigma)$.  We also study the irreducible
components and codimension of $Z(\Sigma) = {\bf A}^n \setminus U(\Sigma)$.

{\em Sections 3 and 4.} We characterize a quasi-smooth hypersurface $X
\subset {\bf P}_\Sigma$ in terms of its defining equation and we
examine relations with $V$-submanifolds and toroidal pairs.  We also
study ${\bf T}(\Sigma)$-linearized sheaves on ${\bf P}_\Sigma$ (where
${\bf T}(\Sigma)$ is the torus acting on ${\bf P}_\Sigma$).

{\em Sections 5, 6 and 7.} The Bott-Steenbrink-Danilov vanishing
theorem for toric varieties was stated without proof in
\cite{dan1,oda}.  We give a proof in the case when ${\bf P}_\Sigma$ is
simplicial, and more generally we prove a vanishing theorem for the
weight filtration on $\Omega^p_{{\bf P}_\Sigma}(\log D)\otimes{\cal
L}$.

{\em Sections 8 and 9.}  Given an ample hypersurface $X
\subset {\bf P}_\Sigma$, we study $H^0({\bf P}_\Sigma,\Omega^p_{{\bf
P}_\Sigma}(X))$.  We use $Cl(\Sigma)$-graded $S(\Sigma)$-modules from
\cite{cox} to describe $\Omega^p_{{\bf P}_\Sigma}$, and we give
explicit generators for the global sections when $p = d$ or $d-1$
(where $d$ is the dimension of ${\bf P}_\Sigma$).

{\em Section 10.}  Given an ample hypersurface $X \subset {\bf
P}_\Sigma$ defined by $f \in S(\Sigma)_\beta$, we show that under the
Hodge filtration on $H^d({\bf P}_\Sigma \setminus X)$, the graded
pieces $Gr_F^pH^d({\bf P}_\Sigma \setminus X)$ are naturally
isomorphic to certain graded pieces of $S(\Sigma)/J(f)$, where $J(f)$
is the Jacobian ideal of $f$.  We then study the primitive cohomology
of $X$ and show how to generalize classical results of Griffiths,
Dolgachev and Steenbrink.

{\em Section 11.}  The first author recently studied the cohomology
of the affine hypersurface $Y = X \cap {\bf T}(\Sigma)$ (see
\cite{bat.var}).  We show how these results can be stated in terms of
the the ideal of $S(\Sigma)$ generated by $z_i \partial f/\partial
z_i$ for $i = 1,\dots,n$.  By looking at the weight filtration on $Y$,
we also get some results on the Hodge components of $X$.

{\em Sections 12 and 13.} The last two sections of the paper give
further results on toric varieties.  First, we generalize the classic
Euler exact sequence and apply it to study $d-1$ forms and the tangent
sheaf on ${\bf P}_\Sigma$.  Then we show how the graded piece
$(S(\Sigma)/J(f))_\beta$ of the Jacobian ring is related to the moduli
of hypersurfaces in ${\bf P}_\Sigma$ defined by $f \in
S(\Sigma)_\beta$.  We use results on the automorphism group of ${\bf
P}_\Sigma$ obtained by the first author in \cite{cox}.

\section{The definition of a simplicial toric variety ${\bf P}_{\Sigma}$}

Let $M$ be a free abelian group of rank $d$, $N = {\rm Hom}(M, {\bf Z})$
the dual group. We denote  by $M_{\bf R}$ (resp. by $N_{\bf R}$) the
${\bf R}$-scalar extension of $M$ (resp. of $N$).

\begin{opr}
{\rm  A convex subset $\sigma \subset N_{\bf R}$ is called a {\em rational
$k$-dimensional simplicial cone} $(k \geq 1)$
if there exist $k$ linearly independent
elements $e_1, \ldots, e_k \in N$ such that
\[ \sigma = \{ \mu_1 e_1 + \cdots + \mu_k e_k \mid
l_i \in {\bf R}, l_i \geq 0 \}. \]
We call $e_1, \ldots, e_k \in N$ {\em integral generators of} $\sigma$
if for every $e_i$ $(1 \leq i \leq k)$ and any non-negative rational
number $\mu$, $\mu \cdot e_i \in N$ only when $\mu \in {\bf Z}$.  The
origin $0 \in N_{\bf R}$ is the {\em rational $0$-dimensional
simplicial cone}, and the set of integral generators of this cone is
empty.}
\end{opr}

\begin{opr}
{\rm A rational simplicial cone $\sigma'$ is called {\em a face} of a
rational simplicial cone $\sigma$ (we write $\sigma' \prec \sigma)$ if the
set of integral generators of $\sigma'$ is a subset of the
set of  integral generators of $\sigma$. }
\label{face}
\end{opr}

\begin{opr}
{\rm A finite set $\Sigma = \{ \sigma_1, \ldots , \sigma_s \}$ of
rational simplicial cones in $N_{\bf R}$ is called
{\em a rational simplicial complete d-dimensional
fan} if the following conditions are satisfied:
\begin{description}
\item[{\rm (i)}] if $\sigma \in  \Sigma$ and $\sigma' \prec \sigma$, then
$\sigma' \in \Sigma$;
\item[{\rm (ii)}] if $\sigma$, $\sigma'$ are in $\Sigma$, then
$\sigma \cap \sigma' \prec \sigma$ and $\sigma \cap \sigma' \prec \sigma'$;
\item[{\rm (iii)}] $N_{\bf R} = \sigma_1 \cup \cdots \cup \sigma_s$.
\end{description}
The set of all $k$-dimensional cones in $\Sigma$ will be denoted by
$\Sigma^{(k)}$.}
\end{opr}

\begin{exam}
{\rm Let $w = \{ w_1, \ldots, w_{d+1}\}$ be the set of
positive integers satisfying the
condition ${\rm gcd}(w_i) =1$. Choose $d+1$ vectors
$ e_1, \ldots, e_{d+1}$ in a
$d$-dimensional real space $V$ such $V$ is spanned by $e_1, \ldots, e_{d+1}$
and there exists the linear relation
\[ w_1 e_1 + \cdots + w_{d+1} e_{d+1} =0. \]
Define $N$ to be the lattice in $V$ consisting of all integral linear
combinations of $e_1, \ldots, e_{d+1}$. Obviously,
$N_{\bf R} = V$. Let $\Sigma(w)$ be the set of
all possible simplicial cones in $V$ generated by proper
subsets of $\{ e_1, \ldots, e_{d+1} \}$. Then $\Sigma(w)$ is an
example of a rational simplicial complete $d$-dimensional fan. }
\label{weig}
\end{exam}

We want to show how every rational simplicial complete $d$-dimensional
fan $\Sigma$ defines a
a compact $d$-dimensional complex algebraic variety ${\bf P}_{\Sigma}$
having only quotient singularities. For instance, if $\Sigma$ is a
fan $\Sigma(w)$ from Example \ref{weig}, then the corresponding variety
${\bf P}_{\Sigma}$ will be the $d$-dimensional weighted projective
space ${\bf P}( w_1, \ldots, w_{d+1})$. The standard defintion  of
${\bf P}( w_1, \ldots, w_{d+1})$  describes it as a quotient of
${\bf C}^{d+1} \setminus \{0\}$ by the diagonal
action of the multiplicative group
${\bf C}^*$:
\[ (z_1, \ldots, z_{d+1} ) \rightarrow
(t^{w_1}z_1, \ldots, t^{w_{d+1}}z_{d+1}),\;\; t \in {\bf C}^*. \]
In Definition \ref{def.tor}, we will show that the toric variety
${\bf P}_{\Sigma}$ can be constructed in a similar manner.  We first
need some definitions.

\begin{opr}
{\rm (\cite{cox}) Let $S(\Sigma) = {\bf C}\lbrack z_1, \ldots, z_n
\rbrack$ be the polynomial ring over ${\bf C}$ with variables $z_1,
\ldots, z_n$, where $\Sigma^{(1)} = \{\rho_1,\dots,\rho_n\}$ are the
$1$-dimensional cones of $\Sigma$.  Then, for $\sigma \in \Sigma$, let
$\widehat{z}_\sigma = \prod_{\rho_i \not\subset \sigma} z_i$, and let
$B(\Sigma) = \langle \widehat{z}_\sigma : \sigma \in \Sigma\rangle \subset
S(\Sigma)$ be the ideal generated by the $\widehat{z}_\sigma$'s. }
\label{def.B}
\end{opr}

\begin{opr}
\label{def.Z}
{\rm Let ${\bf A}^n = {\rm Spec}\, S(\Sigma)$ be the $n$-dimensional
affine space over ${\bf C}$ with coordinates $z_1,\dots,z_n$.  The
ideal $B(\Sigma) \subset S(\Sigma)$ gives the variety
\[ Z(\Sigma) = {\bf V}(B(\Sigma)) \subset {\bf A}^n,\]
and we get the Zariski open set
\[ U(\Sigma) = {\bf A}^n \backslash Z(\Sigma).\]}
\end{opr}

\begin{opr}
{\rm Consider the injective homomorphism $\alpha :
M \rightarrow {\bf Z}^n$ defined by
\[ \alpha(m) = (\langle m,e_1\rangle, \dots, \langle m,e_n\rangle ). \]
The cokernel of this map is ${\bf Z}^n/\alpha(M) \simeq Cl(\Sigma)$,
and we define ${\bf D}(\Sigma) : = {\rm Spec}\,{\bf C} \lbrack Cl(\Sigma)
\rbrack$ to be $(n-d)$-dimensional commutative
affine algebraic $D$-group whose group of characters is isomorphic to
$Cl(\Sigma)$ (see \cite{ham}). }
\label{diag}
\end{opr}

\begin{rem}
{\rm Notice that the finitely generated abelian group
$Cl(\Sigma)$ defines ${\bf D}(\Sigma)$ not only
as an abstract commutative affine algebraic $D$-group, but also
as a canonically  embedded into
$({\bf C}^*)^n = {\rm Spec}\,{\bf C} \lbrack {\bf Z}^n \rbrack$ subgroup
associated
with the surjective homomorphism ${\bf Z}^n \rightarrow
{\bf Z}^n / \alpha(M) \simeq Cl(\Sigma)$.
So we obtain the canonical diagonal action of ${\bf D}(\Sigma)$ on
the affine space ${\bf A}^n$. Obviously, $U(\Sigma)$ is
invariant under this action.}
\end{rem}

Now we are ready to formulate the key theorem.

\begin{theo}
Let $\Sigma$ be a rational simplicial complete $d$-dimensional fan.
Then the canonical action of ${\bf D}(\Sigma)$ on $U(\Sigma)$ has
the geometric quotient
\[ {\bf P}_{\Sigma} = U(\Sigma)/{\bf D}(\Sigma) \]
which is a com\-pact com\-plex $d$-di\-men\-sional al\-geb\-raic
va\-riety having only abel\-ian quotient singularities (and hence
${\bf P}_\Sigma$ is a $V$-manifold).
\label{key.theo}
\end{theo}

\proof  The existence of the quotient $U(\Sigma)/{\bf D}(\Sigma)$ was
proved in the analytic case by Audin \cite{aud} and in the algebraic
case by Cox \cite{cox}.  It remains to show that ${\bf P}_\Sigma$ has
only abelian quotient singularities.  First observe that $Z(\Sigma)$
is defined by the vanishing of $\widehat{z}_\sigma$ for the
$d$-dimensional cones $\sigma \in \Sigma^{(d)}$.  For such a $\sigma$,
set $U_\sigma = \{z \in {\bf A}^n: \widehat{z}_\sigma \ne 0\}$.  Thus
$U(\Sigma) = \bigcup_{\sigma \in \Sigma^{(d)}} U_\sigma$, and it
suffices to show that $U_\sigma/{\bf D}(\Sigma)$ has abelian quotient
singularities.

	If the $1$-dimensional cones of $\Sigma$ are $\Sigma^{(1)} =
\{\rho_1, \ldots, \rho_n \}$, put $G(\Sigma) = \{ e_1, \ldots , e_n
\}$ where $e_i$ is the integral generator of
$\rho_i$.  Now renumber $G(\Sigma)$ so that the
generators of $\sigma$ are $e_1,\dots,e_d$.  Then $U_\sigma = {\bf
C}^d\times({\bf C}^*)^{n-d}$.  To see how ${\bf D}(\Sigma)$ acts on
this set, consider the commutative diagram
\begin{equation}
\begin{array}{ccccccccc}
&&&& 0 && 0 &&\\
&&&& \downarrow && \downarrow &&\\
&&&& {\bf Z}^{n-d} & = & {\bf Z}^{n-d} &&\\
&&&& \downarrow && \downarrow &&\\
0 & \to & M & \to & {\bf Z}^n & \to & Cl(\Sigma) & \to & 0\\
&& \parallel && \downarrow && \downarrow &&\\
0 & \to & M & \to & {\bf Z}^d & \to & Cl(\sigma) & \to & 0\\
&&&& \downarrow && \downarrow &&\\
&&&& 0 && 0 &&
\end{array}
\label{diagram}
\end{equation}
where the map $M \to {\bf Z}^d$ is $m \mapsto (\langle m,e_1\rangle,
\dots, \langle m,e_d\rangle)$.  Since $e_1,\dots,e_d$ are linearly
independent, it follows that $Cl(\sigma)$ is a finite group.  Then the
last column of the diagram gives the exact sequence
\begin{equation}
1 \to {\bf D}(\sigma) \to {\bf D}(\Sigma) \to ({\bf C}^*)^{n-d} \to 1.
\label{eqD}
\end{equation}
Since $U_\sigma = {\bf C}^d\times({\bf C}^*)^{n-d}$ and ${\bf
D}(\sigma)$ acts naturally on ${\bf C}^d$, it follows that the map
${\bf C}^d \to U_\sigma$ defined by $(t_1,\dots,t_d) \mapsto
(t_1,\dots,t_d,1,\dots,1)$ is equivariant.  It is also easy to see
that the induced map ${\bf C}^d/{\bf D}(\sigma) \to U_\sigma/{\bf
D}(\Sigma)$ is an isomorphism.  Since ${\bf D}(\sigma)$ is a finite
abelian group, the theorem is proved. \hfill $\Box$

\begin{rem} {\rm Audin \cite{aud} and Cox \cite{cox} have
shown that ${\bf P}_{\Sigma}$ is isomorphic to the complete toric
variety associated with the simplicial fan $\Sigma$ in the usual sense
of the theory \cite{dan1,oda} (see \cite{cox} for other people who
have discovered this result).  Also, \cite{cox} shows that ${\bf
A}_\sigma = U_\sigma/{\bf D}(\Sigma) \subset {\bf P}_\Sigma$ is the
affine toric open associated to the cone $\sigma \in \Sigma$.  We can
thus use Theorem \ref{key.theo} as the {\em definition} of toric
variety. }
\end{rem}

\begin{opr}
{\rm The complete $d$-dimensional algebraic variety
\[ {\bf P}_{\Sigma} = U(\Sigma)/{\bf D}(\Sigma) \]
is called the {\em toric variety associated with the complete
simplicial fan} $\Sigma$.}
\label{def.tor}
\end{opr}

	This construction of a toric variety makes it easy to see the
torus action.

\begin{opr}
{\rm Denote by ${\bf T}(\Sigma)$ the quotient of $({\bf C}^*)^n$ by
the subgroup ${\bf D}(\Sigma)$. When $\Sigma$ is fixed, we denote
${\bf T}(\Sigma)$ simply by ${\bf T}$. }
\end{opr}

\begin{rem}
{\rm The group ${\bf T}(\Sigma)$ is isomorphic to the torus
$N\otimes{\bf C}^* = ({\bf C}^*)^d$.  Since the group $({\bf C}^*)^n$
is an open subset of $U(\Sigma)$, and it acts canonically on $U(\Sigma)$,
${\bf T}(\Sigma)$ is an open subset of ${\bf P}_{\Sigma}$
having the induced action by regular automorphisms of ${\bf
P}_{\Sigma}$.}
\end{rem}

\section{The structure of $Z(\Sigma) = {\bf A}^n \setminus U(\Sigma)$}

We will next discuss the closed subvariety $Z(\Sigma) = {\bf A}^n
\backslash U(\Sigma)$.  It has an interesting combinatorial
interpretation.

\begin{opr}
{\rm (\cite{bat.class}) Let the $1$-dimensional cones of $\Sigma$ be
$\Sigma^{(1)} = \{\rho_1, \ldots, \rho_n \}$, and put $G(\Sigma) =
\{ e_1, \ldots , e_n \}$ where $e_i$ is the integral generator of
$\rho_i$. We call a subset ${\cal P} =\{ e_{i_1}, \ldots , e_{i_p} \}
\subset G(\Sigma)$ a {\em primitive collection} if $\{e_{i_1}, \ldots
, e_{i_p}\}$ is not the set of generators of a $p$-dimensional
simplicial cone in $\Sigma$, while any proper subset of ${\cal P}$
generates a cone in $\Sigma$.}
\end{opr}

\begin{exam}
{\rm Let $\Sigma$ be a fan $\Sigma(w)$ from Example \ref{weig}. Then
there exists only one primitive collection in $G(\Sigma(w))$ which
coincides with  $G(\Sigma(w))$ itself.}
\label{primweig}
\end{exam}

\begin{opr}
{\rm
Let ${\cal P} =\{ e_{i_1}, \ldots , e_{i_p} \}$ be a primitive collection
in $G(\Sigma)$. Define ${\bf A}({\cal P})$ to be the $(n-p)$-dimensional
affine subspace in ${\bf A}^n$  having the equations
\[ z_{i_1} = \cdots = z_{i_p} = 0. \]}
\end{opr}

\begin{rem}
{\rm Since every primitive collection ${\cal P}$ has at least two elements,
the codimension of ${\bf A}({\cal P})$ is at least $2$.}
\end{rem}

\begin{lem}
Let $Z(\Sigma) \subset {\bf A}^n$ be the variety defined in Definition
\ref{def.Z}.  Then the decomposition of $Z(\Sigma)$ into its
irreducible components is given by
\[ Z(\Sigma) = \bigcup_{\cal P} {\bf A}({\cal P}), \]
where ${\cal P}$ runs over all primitive collections in $G(\Sigma)$.
\label{zdecomp}
\end{lem}

\proof First, it follows from the definition of $Z(\Sigma) = {\bf
V}(\widehat{z}_\sigma : \sigma \in \Sigma)$ that
\begin{equation}
\label{eqz}
Z(\Sigma) = \bigcup_{\cal Q} {\bf A}({\cal Q}),
\end{equation}
where ${\cal Q}$ runs over all subsets ${\cal Q} \subset G(\Sigma)$
which are not the set of generators of any cone in $\Sigma$.  Then
note that the set of all such ${\cal Q}$'s are partially ordered by
inclusion, and ${\cal Q} \subset {\cal Q}'$ implies ${\bf A}({\cal
Q}') \subset {\bf A}({\cal Q})$.  Hence, in the above union, it
suffices to use the minimal ${\cal Q}$'s, which are precisely the
primitive collections.  It follows that these give the irreducible
components of $Z(\Sigma)$. \hfill $\Box$

\begin{rem}
{\rm In \cite{bat.class}, the first author conjectured that for smooth
complete toric varieties, the number of primitive collections could be
bounded in terms of the Picard number $\rho$ ($= n-d$).  Since
primitive collections correspond to irreducible components of
$Z(\Sigma)$ by the above lemma, we can reformulate this conjecture as
follows.}
\end{rem}

\begin{conj} For any $d$-dimensional smooth complete toric variety
with Picard number $\rho$ defined by a complete regular fan $\Sigma$,
there exists a constant $N(\rho)$ depending only on $\rho$ such that
$Z(\Sigma) \subset {\bf A}^n$ has at most $N(\rho)$ irreducible
components.
\end{conj}

	We next study the codimension of $Z(\Sigma)$.  When ${\bf P}_\Sigma$
is a weighted projective space, it follows from Example \ref{primweig}
that $Z(\Sigma) = \{0\}$, which is as small as possible.  It turns out
that in most other cases, $Z(\Sigma)$ is considerably larger.  The
precise result is as follows.

\begin{prop}
\label{prop.codim}
Let ${\bf P}_\Sigma$ be a complete simplicial toric variety of
dimension $d$, and let $Z(\Sigma) \subset {\bf A}^n$ be as above.
Then either
\begin{enumerate}
\item $2 \le {\rm codim}\,Z(\Sigma) \le [{1\over2}d]+1$, or
\item $n = d+1$ and $Z(\Sigma) = \{0\}$.
\end{enumerate}
\end{prop}

\proof To illustrate the range of techniques that can be brought to
bear on this subject, we will give two proofs of this result.  For the
first proof, we make the additional assumption that ${\bf P}_\Sigma$ is
projective.   A projective embedding of ${\bf P}_\Sigma$ is given by a
strictly convex support function on $\Sigma$, which determines
a convex polytope $\Delta \subset M_{\bf R}$ (see \cite{oda}).
Now consider the dual polytope $\Delta^* \subset N_{\bf R}$.
Combinatorially, $\Delta^*$ is closely related to the fan
$\Sigma$---in fact, $\Sigma$ is the cone over $\Delta^*$.

	Let $k = [{1\over2}d]+1$ and assume that ${\rm codim}\,
Z(\Sigma) > k$.  Now pick any subset ${\cal Q} \subset G(\Sigma)$
consisting of $k$ elements (this corresponds to picking $k$ vertices
of $\Delta^*$).  Then ${\rm codim}\, Z(\Sigma) > k$ implies that ${\bf
A}({\cal Q}) \not\subset Z(\Sigma)$.  From equation (\ref{eqz}), it
follows that ${\cal Q}$ must be the set of generators for some face
$\sigma \in \Sigma$.

	In terms of the polytope $\Delta^*$, this shows that every set
of $k$ vertices are the vertices of some face of $\Delta^*$.  Since $k
> [{1\over2}d]$, standard results about convex polytopes (see Chapter
7 of \cite{grunbaum}) imply that $\Delta^*$ is a simplex, which
proves that the number of 1-dimensional cones is $n = d+1$.

	Our second proof, which applies to arbitrary complete
simplicial toric varieties, uses the Stanley-Reisner ring of a certain
monomial ideal associated with $\Sigma$.

\begin{opr}
{\rm Let $\Sigma$ be a complete simplicial fan of dimension $d$. The
{\em Stanley-Reisner ideal} $I(\Sigma)$ is the ideal in $S = S(\Sigma)$
generated by all monomials $z_{i_1} \cdots z_{i_p}$ such that
$\{ e_{i_1}, \ldots, e_{i_p} \}$ is a primitive collection in $G(\Sigma)$.
The quotient $R(\Sigma) = S(\Sigma)/I(\Sigma)$ is called the {\em
Stanley-Reisner ring} of the fan $\Sigma$. }
\end{opr}

As explained in \cite{reisner,stanley}, the ring $R(\Sigma)$ comes
from the simplicial complex given by all subsets of $G(\Sigma)$ which
correspond to cones of $\Sigma$ (this is because $\Sigma$ is a
simplicial fan).  But since $\Sigma$ is also complete, the simplicial
complex is a triangulation of the $d-1$ sphere.  This has some very
strong consequences about the ring $R(\Sigma)$.  For instance, by \S1
of \cite{reisner}, the irreducible components of the affine variety
${\rm Spec}\, R(\Sigma)$ correspond to maximal faces of the simplicial
complex.  It follows that ${\rm Spec}\, R(\Sigma)$ is a union of
$d$-dimensional linear subspaces in ${\bf A}^n$, one for each
$d$-dimensional cone of $\Sigma$.

\begin{rem}
{\rm One should not confuse $Z(\Sigma)$ with ${\rm Spec}\, R(\Sigma)$.
For a weighted projective space as in Example \ref{weig}, $Z(\Sigma) =
\{ 0 \}$, but $I(\Sigma)$ is a principal ideal generated by the
monomial $z_1 \cdots z_{d+1}$, i.e., ${\rm Spec}\, R(\Sigma)$ is the
union of $d+1$ linear subspaces of codimension $1$ in ${\bf A}^{n} =
{\bf A}^{d+1}$. }
\end{rem}

The ring $R(\Sigma)$ has many nice properties.  First, it has the
natural grading by elements of ${\bf Z}_{\geq 0}$ induced from the
${\bf Z}_{\geq 0}$-grading of the polynomials ring $S$.  Second, since
the simplicial complex is a triangulation of the $d-1$ sphere, it
follows from Chapter II, \S5 of \cite{stanley} that $R(\Sigma)$ is a
graded Gorenstein ring of dimension $d$.  Then, by \cite{BE}, the
minimal free resolution
\[ 0 \rightarrow P_{n-d} \stackrel{d_{n-d}}\rightarrow P_{n - d -1}
\stackrel{d_{n-d-1}}\rightarrow \cdots \stackrel{d_{2}}\rightarrow
P_1 \stackrel{d_{0}}\rightarrow P_0 \rightarrow R(\Sigma) \rightarrow 0 \]
of $R(\Sigma)$ as a module over $S$ satisfies the
duality property
\begin{equation}
 P_i \cong {\rm Hom}_S (P_{n-d-i}, P_{n-d})
\label{duality}
\end{equation}
where each $d_i$ is a graded homomorphism of degree $0$ between the
graded $S$-modules $P_i$ and $P_{i-1}$.  Further, the module $P_0$ is
isomorphic to $S$, and since we know the generators of $I(\Sigma)$, we
get an isomorphism
\begin{equation}
P_1 \cong \bigoplus_{{\cal P} \subset G(\Sigma)} S(-|{\cal P}| )
\label{p1}
\end{equation}
where ${\cal P}$ runs over all primitive collections in $G(\Sigma)$,
and $|{\cal P}|$ is the cardinality of ${\cal P}$.  Finally, the
duality (\ref{duality}) shows that $P_{n-d}$ is a free $S$-module of
rank 1.

	We claim that $P_{n-d} \cong S(-n)$.  For each $i$ between $0$
and $n-d$, let $h(i)$ be the minimal integer $h$ such that $P_i$ has a
nonzero element of degree $h$.  The minimality of the free resolution
implies that $0 = h(0) < h(1) < \cdots < h(n-d-1) < h(n-d)$.  Since
$P_{n-d}$ has rank 1, we have $P_{n-d} = S(-h(n-d))$.  Hence it
suffices to show $h(n-d) = n$.  The Hilbert-Poincare series of $S(-j)$
is $H(S(-j),t) = t^j/(1-t)^n$, so that the free resolution of
$R(\Sigma)$ implies that
\[H(R(\Sigma),t) = \sum_{i=0}^{n-d} (-1)^i H(P_i,t) =
{\hbox{polynomial of degree $h(n-d)$} \over (1-t)^n}.\]
However, Theorem 1.4 of Chapter II of \cite{stanley} shows that
\[ H(R(\Sigma),t) = {\hbox{polynomial of degree $d$} \over (1-t)^d}.\]
Comparing these two expressions, we conclude $h(n-d) = n$.

	Once we know $P_{n-d} = S(-n)$, the isomorphism (\ref{p1}) and
the duality (\ref{duality}) give an isomorphism
\begin{equation}
P_{n-d-1} \cong \bigoplus_{{\cal P} \subset G(\Sigma)} S(-n+|{\cal P}|).
\label{pn-d-1}
\end{equation}
Assume now that $n > d+ 1$. Then $h(n-d-1) > \cdots > h(1)$
implies
\begin{equation}
h(n -d -1) \geq h(1) + n - d - 2.
\label{ineq}
\end{equation}
{}From the isomorphisms (\ref{p1}) and (\ref{pn-d-1}), we see that
there exist primitive collections ${\cal P}_1$ and ${\cal P}_2$ such
that
\[ h(n-d-1) =  n - |{\cal P}_1|\quad\hbox{and}\quad h(1) = |{\cal P}_2|. \]
Then the inequality (\ref{ineq}) implies that
\[  |{\cal P}_1|+ |{\cal P}_2| \leq  d + 2, \]
and thus
\[ \min(|{\cal P}_1|,|{\cal P}_2|) \le [\textstyle{1\over2}d] + 1.\]
On the other hand, Lemma \ref{zdecomp} implies
\[ \min_{{\cal P} \subset G(\Sigma)} |{\cal P}| = \min_{{\cal P}
\subset G(\Sigma)} {\rm codim}\, {\bf A}({\cal P}) = {\rm codim}\,
Z(\Sigma). \]
This proves ${\rm codim}\, Z(\Sigma) \le [{1\over2}d] + 1$ when $n >
d+1$.

On the other hand, the equality $n = d+1$ is possible only if the
minimal projective resoluton of $R(\Sigma)$ consists of $P_0 \cong S$
and $P_1 \cong S(-n)$, which means that $I(\Sigma)$ is the principal
ideal generated by the monomial $z_1 \cdots z_n$ of degree $n$.  This
implies that $\{z_1,\dots,z_n\}$ is the unique primitive collection,
and we obtain $Z(\Sigma) = \{0 \}$.  This completes the proof of the
proposition. \hfill$\Box$
\medskip

	The condition $n = d+1$ is closely related to ${\bf P}_\Sigma$
being a weighted projective space.

\begin{lem} Let $\Sigma$ be a complete simplicial fan with $n = d+1$
1-dimensional cones.  Then there is a weighted projective space ${\bf
P}(w_1,\dots,w_{d+1})$ and a finite surjective morphism
\[ {\bf P}(w_1,\dots,w_{d+1}) \to {\bf P}_\Sigma. \]
Furthermore, if $G(\Sigma) = \{e_1,\dots,e_{d+1}\}$, then the
following are equivalent:
\begin{enumerate}
\item ${\bf P}_\Sigma$ is a weighted projective space.
\item ${\bf D}(\Sigma) \cong {\bf C}^*$ (or equivalently, $Cl(\Sigma)
\cong {\bf Z}$).
\item $e_1,\dots,e_{d+1}$ generate $N$ as a ${\bf Z}$-module.
\end{enumerate}
\label{lem.nd}
\end{lem}

\proof First note that (1) $\Rightarrow$ (2) is immediate.  Further,
if (3) holds, then Example \ref{weig} shows that ${\bf P}_\Sigma$ is a
weighted projective space.  It remains to show (2) $\Rightarrow$ (3).
But if $Cl(\Sigma) \cong {\bf Z}$, then taking the dual of the exact
sequence
\[ 0 \to M \stackrel{\alpha}\rightarrow {\bf Z}^{d+1} \to {\bf Z}
\to 0\]
from Definition \ref{diag} gives an exact sequence
\[ 0 \to {\bf Z} \to {\bf Z}^{d+1} \to N \to 0,\]
and it follows that $e_1,\dots,e_{d+1}$ generate $N$.

	Finally, given an arbitrary fan $\Sigma$ with $n = d+1$, let
$N' \subset N$ be the sublattice generated by $e_1,\dots,e_{d+1}$.
Then the fan $\Sigma$ induces a fan $\Sigma'$ for the lattice $N'$,
and the natural inclusion $(N',\Sigma') \to (N,\Sigma)$ induces a
finite surjection ${\bf P}_{\Sigma'} \to {\bf P}_\Sigma$ by Corollary
1.16 of \cite{oda}.  By the above, ${\bf P}_{\Sigma'}$ is a weighted
projective space, and the proposition is proved. \hfill $\Box$
\medskip

	Lemma \ref{lem.nd} implies that a smooth complete toric
variety with $n=d+1$ is ${\bf P}^d$ (this fact is well-known).  Thus
we get the following corollary of Proposition \ref{prop.codim}.

\begin{coro} Let ${\bf P}_\Sigma$ be a smooth complete toric variety.
Then either
\begin{enumerate}
\item $2 \le {\rm codim}\, Z(\Sigma) \le [{1\over2}d]+1$, or
\item ${\bf P}_\Sigma = {\bf P}^d$.
\end{enumerate}
\end{coro}

We end this section with some comments about combinatorially
equivalent fans.

\begin{opr}
{\rm Two rational simplicial complete $d$-dimensional fans $\Sigma$ and
$\Sigma'$ are called {\em combinatorially equivalent} if there exists a
bijective mapping $\Sigma \rightarrow \Sigma'$ respecting the
face-relation ``$\prec$" (see Definition \ref{face}). }
\end{opr}

\begin{rem}
{\rm It is easy to see that the closed subset $Z(\Sigma) \subset {\bf
A}^n$ depends only on the combinatorial structure of $\Sigma$, i.e.,
for two combinatorially equivalent fans $\Sigma$ and $\Sigma'$, we
can assume that we are in the same affine space ${\bf A}^n$, and then
we have $Z(\Sigma) = Z(\Sigma')$ and hence $U(\Sigma) = U(\Sigma')$. }
\end{rem}

It follows that all toric varieties coming from combinatorially
equivalent fans are quotients of the same open set $U \subset {\bf
A}^n$, though the group actions will differ.  For example, all
$d$-dimensional weighted projective spaces come from combinatorially
equivalent fans and are quotients of $U = {\bf A}^{d+1} - \{0\}$.

\section{Quasi-smooth hypersurfaces}

	Throughout this section, let ${\bf P} = {\bf P}_\Sigma$ be a
fixed $d$-dimensional complete simplicial toric variety.  The action
of ${\bf D} = {\bf D}(\Sigma)$ on ${\bf A}^n$ induces an action on $S
= S(\Sigma) = {\bf C}[z_1,\dots,z_n]$.  The decomposition of this
representation gives a grading on $S$ by the character group
$Cl(\Sigma)$ of ${\bf D}$.  A polynomial $f$ in the graded piece of
$S$ corresponding to $\beta \in Cl(\Sigma)$ is said to be ${\bf
D}$-homogeneous of degree $\beta$.

	 Such a polynomial $f$ has a zero set ${\bf V}(f) \subset {\bf
A}^n$, and ${\bf V}(f)\cap U(\Sigma)$ is stable under ${\bf D}$ and
hence descends to a hypersurface $X \subset {\bf P}$ (this is because
${\bf P}$ is a geometric quotient---see \cite{cox} for more details).
We call ${\bf V}(f) \subset {\bf A}^n$ the {\em affine quasi-cone} of
$X$.

\begin{opr}
{\rm If a hypersurface $X \subset {\bf P}$ is defined by a ${\bf
D}$-homogeneous polynomial $f$, then we say that $X$ is {\em
quasi-smooth} if the affine quasi-cone ${\bf V}(f)$ is smooth outside
$Z = Z(\Sigma) \subset {\bf A}^n$.  }
\label{quasi}
\end{opr}

	To see what this definition says about the hypersurface $X$,
we need the concept of a $V$-submanifold of a $V$-manifold.  As usual,
a $d$-dimensional variety $W$ is a $V$-manifold if for every point $p
\in W'$, there is an analytic isomorphism of germs $({\bf C}^d/G,0)
\cong (W,p)$ where $G \subset GL(d,{\bf C})$ is a finite small
subgroup.  (Recall that being {\em small} means there are no elements
with 1 as an eigenvalue of multiplicity $d-1$, i.e., no complex
rotations other than the identity).  In this case, we say that $({\bf
C}^d/G,0)$ is a {\em local model} of $W$ at $p$.

\begin{opr}
{\rm If a $d$-dimensional variety $W$ is a $V$-manifold, then a
subvariety $W' \subset W$ is a {\em $V$-submanifold} if for every point
$p \in W'$, there is a local model $({\bf C}^d/G,0) \cong (W,p)$ such
that $G \subset GL(d,{\bf C})$ is a finite small subgroup and the
inverse image of $W'$ in ${\bf C}^d$ is smooth at $0$.  }
\end{opr}

\begin{rem}
{\rm In \cite{prill}, Prill proved that $G$ is determined up to
conjugacy by the germ $(W,p)$.  Hence, if $W' \subset W$ is a
$V$-submanifold for one local model at $p$, then it is a
$V$-submanifold for all local models at $p$.  It is also easy to see
that a $V$-submanifold of a $V$-manifold $W$ is again a $V$-manifold.
However, the converse is false: a subvariety of $W$ which is a
$V$-manifold need not be a $V$-submanifold.  This is because the
singularities of a $V$-submanifold are intimately related to the
singularities of the ambient space.}
\end{rem}

	We will omit the proof of the following easy lemma.

\begin{lem} If $G \subset GL(d,{\bf C})$ is a small and finite, then a
subvariety $W' \subset {\bf C}^d/G$ is a $V$-submanifold if and only if
the inverse image of $W'$ in ${\bf C}^d$ is smooth.
\label{local}
\end{lem}

	We can now describe when a hypersurface of the toric variety
${\bf P}$  is a $V$-submanifold.

\begin{prop} If a hypersurface $X \subset {\bf P}$ is defined by
a ${\bf D}$-homogeneous polynomial $f$, then $X$ is quasi-smooth if
and only if $X$ is a $V$-submanifold of ${\bf P}$.
\label{prop.quasi}
\end{prop}

\proof We will use the notation of the proof of Theorem
\ref{key.theo}.  Thus, let $\sigma$ be a $d$-dimensional cone of
$\Sigma$ and assume that $\sigma$ is generated by $e_1,\dots,e_d$.
As we saw in the proof of Theorem~\ref{key.theo}, $\sigma$ gives an
affine open set ${\bf A}_\sigma = {\bf C}^d/{\bf D}(\sigma) \cong
U_\sigma/{\bf D}$ of ${\bf P}$.

	We claim that ${\bf D}(\sigma) \subset ({\bf C}^*)^d$ is a
small subgroup.  Suppose that $g = (\lambda,1,\dots,1) \in {\bf
D}(\sigma)$.  We can regard $g$ as a homomorphism $g : Cl(\sigma) \to
{\bf C}^*$, and diagram (\ref{diagram}) shows that $g = 1$ on the
image of $(\langle m,e_1\rangle,\dots,\langle m,e_d\rangle)$ whenever
$m \in M$.  This means that $\lambda^{\langle m,e_1\rangle} = 1$ for
all $m \in M$.  Since $e_1$ is not the multiple of any element of $N$
and $M = {\rm Hom}(N,{\bf Z})$, it follows that $\lambda = 1$ and
hence $g$ is the identity.

	Since ${\bf A}_\sigma = {\bf C}^d/{\bf D}(\sigma)$, Lemma
\ref{local} implies that $X\cap {\bf A}_\sigma$ is a $V$-submanifold
if any only if the inverse image of $X\cap{\bf A}_\sigma$ in ${\bf
C}^d$ is smooth.  Since the map ${\bf C}^d/{\bf D}(\sigma) \cong
U_\sigma/{\bf D}$ is induced by $(t_1,\dots,t_d) \mapsto
(t_1,\dots,t_d,1,\dots,1)$, it follows that the inverse image of
$X\cap{\bf A}_\sigma$ is ${\bf C}^d$ is defined by $g = 0$, where
\[ g(t_1,\dots,t_d) = f(t_1,\dots,t_d,1,\dots,1).\]
Thus $X\cap{\bf A}_\sigma$ is a $V$-submanifold if and only if the
subvariety $g = 0$ is smooth in ${\bf C}^d$.

	We need to relate this to the smoothness of ${\bf V}(f)\cap
U_\sigma$.  Recall the decomposition $U_\sigma = {\bf C}^d\times({\bf
C}^*)^{n-d}$.  Using the surjection ${\bf D} \to ({\bf C}^*)^{n-d}$
from equation (\ref{eqD}), we see that ${\bf V}(f)\cap U_\sigma$ is
smooth if and only if it is smooth at all points of the form
$(t_1,\dots,t_d,1,\dots,1)$.  Now comes the key observation.

\begin{lem} ${\bf V}(f)\cap U_\sigma$ is smooth at $t =
(t_1,\dots,t_d,1,\dots,1)$ if and only if one of the partials
$f_{z_i}(t)$ is nonzero for some $1 \le i \le d$.
\label{key.lem}
\end{lem}

\begin{rem} {\rm This lemma and the previous paragraphs show that
${\bf V}(f)\cap U_\sigma$ is smooth if and only if $X\cap{\bf
A}_\sigma$ is a $V$-submanifold.  Since the quasi-cone of $X$ is
smooth outside $Z$ if and only if ${\bf V}(f)\cap U_\sigma$ is smooth
for all $d$-dimensional $\sigma$, Proposition \ref{prop.quasi} follows
from Lemma \ref{key.lem}.}
\end{rem}

\proof  Suppose $f_{z_i}(t) = 0$ for $1 \le i \le d$.  Now take any
$j > d$.  Since $e_1,\dots,e_d$ are a basis of $N_{\bf Q}$, we can
write $e_j = \sum_{i=0}^d \phi_i e_i$.  By Lemma \ref{euler.lem}
below, there is a constant $\phi(\beta)$ such that
\[ \phi(\beta) f = z_j{\partial f \over \partial z_j} - \sum_{i=1}^d
\phi_i z_i {\partial f \over \partial z_i}.\]
Evaluating this at $t \in {\bf V}(f)$ and using $f_{z_i}(t) = 0$ for
$1 \le i \le d$, we see that $f_{z_j}(t) = 0$ for all $1 \le j \le
n$.  Hence, in order to be smooth at $t$, at least one of the first
$d$ partials must be nonvanishing. \hfill$\Box$
\medskip

	To complete the proof of Proposition \ref{prop.quasi}, we need
to prove the following lemma.

\begin{lem}
\label{euler.lem}
Suppose that we have complex numbers $\phi_1,\dots,\phi_n$
with the property that $\sum_{i=0}^n \phi_i e_i = 0$ in $N_{\bf C}$.
Then, for any class $\beta \in Cl(\Sigma)$, there is a constant
$\phi(\beta)$ with the property that for any ${\bf D}$-homogeneous
polynomial $f \in S$ of degree $\beta$, we have
\[ \phi(\beta) f = \sum_{i=1}^n \phi_i z_i {\partial f \over \partial
z_i}.\]
\end{lem}

\proof From $\phi_1,\dots,\phi_n$, we get a map $\tilde\phi : {\bf
Z}^n \to {\bf C}$ defined by $(a_1,\dots,a_n) \mapsto \sum_{i=1}^n \phi_i
a_i$.  Furthermore, under the map $\alpha: M \to {\bf Z}^n$, note that
$m \in M$ maps to $\sum_{i=1}^n \phi_i \langle m,e_i\rangle =
\big\langle m,\sum_{i=1}^n \phi_i e_i\big\rangle = 0$.  Thus
$\tilde\phi$ induces a map $\phi : {\bf Z}^n/\alpha(M) \cong
Cl(\Sigma) \to {\bf C}$.

	Note that every monomial $f = \prod_{i=1}^n z_i^{a_i}$ of $S$
is ${\bf D}$-homogeneous.  If we let $\beta = \deg f = \sum_{i=1}^n
a_i \deg z_i$, then the desired formula for $\phi(\beta) f$ follows
immediately.  By linearity, the formula then holds for all ${\bf
D}$-homogeneous poynomials of degree $\beta$. \hfill$\Box$
\medskip

	In light of this lemma, we make the following definition.

\begin{opr} {\rm The identity $\phi(\beta) f = \sum_{i=1}^n \phi_i z_i
f_{z_i}$ from Lemma \ref{euler.lem} is called the {\em Euler formula}
determined by $\phi_1,\dots,\phi_n$.}
\label{euler.def}
\end{opr}

\begin{rem} {\rm It follows easily from the proof of Lemma
\ref{euler.lem} that the set of all Euler formulas form the vector
space ${\rm Hom}(Cl(\Sigma), {\bf C})$.  For example, if ${\bf P}$ is
a weighted projective space, then ${\rm Hom}(Cl(\Sigma), {\bf C})
\cong {\bf C}$, where a basis is given by the usual Euler formula for
homogeneous polynomials.  Note also that ${\rm Hom}(Cl(\Sigma),{\bf
C}) = {\rm Lie}({\bf D})$, so that an Euler formula can be regarded as
a vector field $\sum_{i=1}^n \phi_i z_i{\partial\over\partial z_i}$
which is tangent to the orbits of ${\bf D}$.}
\end{rem}

	We will end this section with a discussion of how quasi-smooth
hypersurfaces of ${\bf P}$ relate to Danilov's theory of ``toroidal
pairs''
\cite{dan3}.

\begin{opr} {\rm Let $W'$ be a subvariety of an algebraic variety $W$.
Then the pair $(W,W')$ is {\em simplicially toroidal} if for every
point $p \in W$, there is a simplicial cone $\sigma$ such that the
germ of $(W,W')$ at $p$ is analytically isomorphic to the germ of
$({\bf A}_\sigma,D)$ at the origin, where ${\bf A}_\sigma$ is the
toric variety of $\sigma$ and $D$ is an irreducible torus-invariant
subvariety of ${\bf A}_\sigma$.  }
\label{def.dan}
\end{opr}

	For hypersurfaces of ${\bf P}$, this concept is equivalent to
being quasi-smooth.

\begin{prop} A hypersurface $X \subset {\bf P}$ is quasi-smooth if and
only if the pair $({\bf P},X)$ is simplicially toroidal.
\label{toroid.qsmooth}
\end{prop}

\proof First suppose we have a pair $({\bf A}_\sigma,D)$ as in
Definition \ref{def.dan}.  We know that ${\bf A}_\sigma \cong {\bf
C}^d/{\bf D}(\sigma)$, and since the inverse image of $D$ is a
coordinate subspace, we see that $D$ is a $V$-submanifold of ${\bf
A}_\sigma$.  Thus $({\bf P},X)$ simplicially toroidal implies that $X$
is a $V$-submanifold and hence quasi-smooth by Proposition
\ref{prop.quasi}.

	Conversely, if $X$ is quasi-smooth, then it is a
$V$-submanifold of ${\bf P}$.  Thus, given $p \in X$, there is a local
model $({\bf C}^d/G,0) \cong ({\bf P},p)$ such that $G \subset ({\bf
C}^*)^d$ is small and finite, and the inverse image $Y \subset {\bf
C}^d$ of $X \subset {\bf P}$ is smooth at the origin.  Then $Y$ is
defined by some equation $h = 0$, and since $Y$ is $G$-invariant,
there is a character $\lambda : G \to {\bf C}^*$ such that $h(g\cdot
t) = \lambda(g)h(t)$ for all $g \in G$.  Since $Y$ is smooth at the
origin, we can also assume that the partial derivative $h_{z_1}(0)$ is
nonvanishing.

	Then the map $\phi : {\bf C}^d \to {\bf C}^d$ defined by
$(t_1,\dots,t_d) \mapsto (h(t_1,\dots,t_d),t_2,\dots,t_d)$ is a local
analytic isomorphism which carries $Y$ to a coordinate hyperplane.
Note that $\phi$ is equivariant provided that $g = (g_1,\dots,g_d) \in
G$ acts on the target space via
\[ g\cdot (t_1,\dots,t_d) = (\lambda(g)t_1,g_2 t_2,\dots,g_d t_d).\]
This gives a map $\psi : G \to ({\bf C}^*)^d$.

	We thus have a local model $({\bf C}^d/\psi(G),0)$ where the
inverse image of $X$ is a coordinate hyperplane.  Since the torus
$({\bf C}^*)^d/\psi(G)$ acts on the affine variety ${\bf
C}^d/\psi(G)$, it follows from Theorem 1.5 of \cite{oda} that ${\bf
C}^d/\psi(G)$ is an affine toric variety.  Thus ${\bf C}^d/\psi(G)
\cong {\bf A}_\sigma$ for some cone $\sigma$, and note that $\sigma$
must be simplicial.  Finally, the image in ${\bf C}^d/\psi(G)$ of a
coordinate hyperplane is an irreducible torus-invariant subvariety.
This proves that $({\bf P},X)$ has the appropriate simplicial toroidal
local model, and the proposition is proved.
\hfill $\Box$

\begin{rem}
{\rm Let ${\bf A}_{\sigma} = {\rm Spec}\, {\bf C} \lbrack
\check{\sigma} \cap M \rbrack$ be an affine $d$-dimensional toric
variety correponding to a $d$-dimensional rational simplicial cone
$\sigma \in \Sigma$, where $\check{\sigma}$ is the dual cone in $M_{\bf
R}$.  There are two important cases where one can explicitly describe
local toroidal models for a quasi-smooth hypersurface $X$ in ${\bf
A}_{\sigma}$:

{\sc Case I.} $X$ has transversal intersections with all orbits ${\bf
T}_\tau \subset {\bf A}_{\sigma}$ of the action of the torus ${\bf T}$
on ${\bf A}_{\sigma}$. Then at the point of intersection of $X$ with a
$1$-dimensional stratum ${\bf T}_{\tau_0}$ corresponding to a
$(d-1)$-dimensional face $\tau_0 \prec \sigma$, the local toroidal
model of $X$ is the $(d-1)$-dimensional affine toric variety ${\bf
A}_{\tau_0}$.

{\sc Case II.} $X$ contains the single closed ${\bf T}_\sigma$-orbit
$p_{\sigma} \in {\bf A}_{\sigma}$ and is ``tangent" to a closure of a
$(d-1)$-dimensional ${\bf T}_\sigma$-orbit corresponding to a
$1$-dimensional cone $\rho \prec \sigma$. In this case, the local
toroidal model of $X$ at $p_{\sigma}$ is the $(d-1)$-dimensional
affine toric variety ${\bf A}_{\sigma/\rho}$, where $\sigma/\rho$ is
the $(d-1)$-dimensional projection of the cone $\sigma \subset N_{\bf
R}$ to the quotient $N_{\bf R} / {\bf R}\rho$.  }
\end{rem}

\section{${\bf T}$-linearized sheaves on  ${\bf P}$}

In this section we will assume that ${\bf P}$ is a complete toric
variety.  Recall from \cite{mam} the notion of a linearization of a
sheaf on an algebraic variety having a regular action of an algebraic
group.  We will apply this to sheaves on ${\bf P}$ with its action by
the torus ${\bf T}$.

\begin{opr}
{\rm Let $\mu_{\bf t} : {\bf P} \rightarrow {\bf P}$ be the
automorphism of ${\bf P}$ defined by $t \in {\bf T}$. Then a {\em {\bf
T}-linearization of a sheaf ${\cal E}$} is a family of
isomorphisms
\[ \phi_t\; : \; \mu_{t}^*\, {\cal E} \cong {\cal E} \]
satisfying the co-cycle condition
\[ \phi_{{t_1}\cdot {t}_2} = \phi_{{t}_2} \circ \mu^*_{{t}_2}
\phi_{{t}_1},\;{\rm for\; all} \;{t}_1, { t}_2 \in {\bf T}. \]}
\end{opr}

\begin{rem}
{\rm If ${\cal E}$ is a ${\bf T}$-linearized sheaf on ${\bf P}$, then
for any ${\bf T}$-invariant open subset ${U} \subset {\bf P}$, the
group ${\bf T}$ has the natural linear representation in the space of
global sections $H^0({ U}, {\cal E})$. Thus $H^0(U, {\cal E})$
splits into a direct sum of subspaces $H^0(U, {\cal E})_{m}$
corresponding to characters $ m \in M$ of ${\bf T}$. }
\end{rem}

\begin{opr}
\label{def.polytope}
{\rm (\cite{kempf})
Let ${\cal E}$ be a ${\bf T}$-linearized sheaf on ${\bf P}$. Then the
polytope
\[ \Delta({\cal E}) = {\rm Conv}\, \{ m \in M : H^0({\bf P}, {\cal E})_{m}
\neq 0 \} \]
is called the {\em support polytope for} ${\cal E}$.}
\end{opr}

\begin{opr}
{\rm Let ${\cal L}$ be a ${\bf T}$-linearized invertible sheaf on
${\bf P}$, $\sigma \in \Sigma^{(d)}$ a $d$-dimensional cone in
$\Sigma$. We denote by $m_{\sigma}({\cal L})$ the unique element of
$M$ with the property that
\[ \{m \in M : H^0(X_{\sigma}, {\cal L})_m \ne 0\} = m_\sigma({\cal
L}) + \check\sigma\cap M,\]
where ${\bf A}_\sigma \subset {\bf P}$ is the affine toric variety
determined by $\sigma$. }
\end{opr}

\begin{rem}
{\rm The existence of $m_\sigma({\cal L})$ follows from \S6.2 of
\cite{dan1}, and it is unique because $\sigma$ is $d$-dimensional.
Note also that the mapping $h : N_{\bf R} \to {\bf R}$ defined by
$h(n) = \langle m_\sigma({\cal L}),n\rangle$ (for $n \in \sigma$) is the
support function for the invertible sheaf ${\cal L}$.

We should also mention that every invertible sheaf ${\cal L}$ on ${\bf
P}$ has a ${\bf T}$-linearization and that any two linearizations of
${\cal L}$ differ by a homomorphism ${\bf T} \to {\bf C}^*$, i.e., by
an element of $M$.  }
\end{rem}

\begin{prop}
\label{delta}
{\rm (\cite{dan1})}
Let ${\cal L}$ be a ${\bf T}$-linearized invertible sheaf on ${\bf
P}$. Then the corresponding polyhedron $\Delta = \Delta({\cal L})$ is
equal to the intersection
\[ \bigcap_{\sigma \in \Sigma^{(d)}} ( m_{\sigma}({\cal L}) +
\check\sigma). \]
\end{prop}

\begin{opr}
{\rm Let $\tau$ be a $k$-dimensional cone in $\Sigma$, and let ${\rm
St}\,(\tau)$ the set of all cones $\sigma \in \Sigma$ such that $\tau
\prec \sigma$. Consider the $(d-k)$-dimensional fan $\Sigma(\tau)$
consisting of projections of cones $\sigma \in {\rm St}\,(\tau)$ into
$N_{\bf R} / {\bf R}\tau$. The corresponding $(d-k)$-dimensional
toric subvariety in ${\bf P}$ will be denoted by ${\bf P}_{\tau}$.}
\label{def.orbit}
\end{opr}

\begin{rem}
{\rm The toric subvariety ${\bf P}_{\tau}$ is the closure of a
$(d-k)$-dimensional orbit ${\bf T}_{\tau}$ of ${\bf T}$. Any ${\bf
T}_{\tau}$-linearized sheaf ${\cal E}$ on ${\bf P}_{\tau}$ can be
considered as a ${\bf T}$-linearized sheaf on ${\bf P}$. }
\end{rem}

\begin{opr}
{\rm Let ${\cal L}$ be a ${\bf T}$-linearized ample invertible sheaf
on ${\bf P}$.  For any $k$-dimensional cone in $\tau \in \Sigma$, we
denote by $\Delta_{\tau}$ the face of $\Delta$ of codimension $k$
defined as
\[ \bigcap_{\sigma \in \Sigma^{(d)},\, \tau \prec \sigma}
( m_{\sigma}({\cal L}) + \check\sigma \cap \tau^{\perp}). \]}
\end{opr}

The next statement follows immediately from the ampleness criterion for
invertible sheaves \cite{dan1}.

\begin{prop}
One has the one-to-one correspondence between $(d-k)$-dimensional
faces $\Delta_{\tau}$
of the polytope $\Delta = \Delta({\cal L})$ and $k$-dimensional
cones $\tau \in \Sigma$ reversing the face-relation. Moreover,
the $(d-k)$-dimensional polytope $\Delta_{\tau}$ is the support polytope
for  the ${\bf T}$-linearized sheaf ${\cal O}_{{\bf P}_{\tau}}
\otimes {\cal L}$.
\label{ample}
\end{prop}

	We next study the relation between $H^0({\bf P},{\cal L})$ and
the coordinate ring $S = {\bf C}[z_1,\dots,z_n]$.

\begin{lem} If ${\cal L}$ is a ${\bf T}$-linearized invertible sheaf
on ${\bf P}$, and let $\beta \in Cl(\Sigma)$ be the class of ${\cal
L}$.  Then there is a natural isomorphism
\[ H^0({\bf P},{\cal L}) \cong S_\beta,\]
where $S_\beta$ is the graded piece of $S$ corresponding to $\beta$.
This isomorphism is determined uniquely up to a nonzero constant in
${\bf C}$.
\label{lem.iso}
\end{lem}

\proof Since ${\bf P}$ is complete, there is a one-to-one
correspondance between ${\bf T}$-linearized invertible sheaves and
${\bf T}$-invariant Cartier divisors (see \S2.2 of \cite{oda}).  Thus
there is a ${\bf T}$-invariant Cartier divisor such that ${\cal L}
\cong {\cal O}_{\bf P}(D)$ as ${\bf T}$-linearized sheaves.  Note that
this isomorphism is unique up to a nonzero constant.  However, in
\cite{cox}, it is shown that $D$ determines an isomorphism $H^0({\bf
P},{\cal O}_{\bf P}(D)) \cong S_\beta$, and the lemma follows
immediately.\hfill$\Box$

\begin{rem}
{\rm If in addition ${\bf P}$ is simplicial, then the polynomial $f
\in S_\beta$ corresponding to a global section of ${\cal L}$
determines a hypersurface $X \subset {\bf P}$ as in \S3, and one can
check that $X$ is exactly the zero section of the global section.}
\end{rem}

\begin{opr}
\label{def.nondeg}
{\rm Let $f$ be a global section of an ample invertible sheaf ${\cal
L}$ on ${\bf P}$. Then the hypersurface $X = \{ p \in {\bf P} : f(p) =
0 \}$ is called {\em nondegenerate} if for any $\tau \in \Sigma$,
the affine hypersurface $X \cap {\bf T}_{\tau}$ is a smooth subvariety
of codimension $1$ in ${\bf T}_{\tau}$. }
\end{opr}

\begin{rem}
{\rm The nondegeneracy of a global section $f$ is equivalent to $f$
being $\Delta$-regular, as defined in \cite{bat.var} (where $\Delta =
\Delta({\cal L})$).  For a proof of this, see \cite{khov}.}
\end{rem}

\begin{prop}
Let $f$ be a generic global section of an ample ${\bf T}$-linearized
invertible sheaf ${\cal L}$ on ${\bf P}$. Then $X = \{p \in {\bf P} :
f(p) = 0 \}$ is a nondegenerate hypersurface.  Moreover, every
nondegenerate hypersurface $X \subset {\bf P}$ is quasi-smooth.
\label{prop.generic}
\end{prop}

\proof  As observed in \cite{dan2}, the first part of the statement
follows from Bertini's theorem, and it was proved in \cite{dan1} that
$X$ is simplicially toroidal. Thus, it remains to apply Proposition
\ref{toroid.qsmooth}. \hfill $\Box$

\begin{rem}
{\rm One should remark that a quasi-smooth hypersurface in a toric
variety ${\bf P}$ need not be nondegenerate.}
\end{rem}

	We will conclude this section by studying the relation between
${\bf T}$-linearized sheaves on ${\bf P}$ and graded $S$-modules.  It
is known (see \cite{cox}) that every $Cl(\Sigma)$-graded $S$-module
$F$ determines a quasi-coherent sheaf $\widetilde{F}$ on ${\bf P}$.
What extra structure on $F$ is needed in order to induce a ${\bf
T}$-linearization on $\widetilde{F}$?  To state the answer, note that
$S$ has a natural grading by ${\bf Z}^n$ which is compatible with its
grading by $Cl(\Sigma)$ via the map ${\bf Z}^n \to Cl(\Sigma)$ from
Definition \ref{diag}.  Then any ${\bf Z}^n$-graded module can be
regarded as a $Cl(\Sigma)$-graded module.

\begin{prop}
If $F$ is a ${\bf Z}^n$-graded $S$-module, then the sheaf
$\widetilde{F}$ on ${\bf P}$ has a natural ${\bf T}$-linearization.
Furthermore, if ${\bf P}$ is simplicial and $N$ is generated by
$e_1,\dots,e_n$, then every ${\bf T}$-linearized quasi-coherernt sheaf
on ${\bf P}$ arises in this way.
\end{prop}

\proof Given $\sigma \in \Sigma$, let $S_\sigma$ be the localization
of $S$ at $\widehat{z}_\sigma = \prod_{e_i \notin \sigma} z_i$, and
let $(S_\sigma)_0$ be the elements of degree $0$ with respect to
$Cl(\Sigma)$.  Then the module $(F\otimes S_\sigma)_0$ determines a
sheaf on the affine piece ${\bf A}_\sigma = {\rm Spec}((S_\sigma)_0)$
of ${\bf P}$, and, as explained in \S3 of \cite{cox}, these sheaves
patch to give $\widetilde{F}$.

	Since $S_\sigma$ and $F\otimes S_\sigma$ have natural ${\bf
Z}^n$ gradings, the groups $(S_\sigma)_0$ and $(F\otimes S_\sigma)_0$
have natural gradings by $M$ (which is the kernel of ${\bf Z}^n \to
Cl(\Sigma)$).  The $M$ grading on $(S_\sigma)_0$ determines the action
of ${\bf T}$ on ${\bf A}_\sigma$, and the $M$ grading on $(F\otimes
S_\sigma)_0$ then determines a ${\bf T}$-linearization.  These clearly
patch to give a ${\bf T}$-linearization of $\widetilde{F}$.

	To prove the final part of the proposition, note that
$Cl(\Sigma)$ is torsion free since $e_1,\dots,e_n$ generate $N$.  Thus
the map ${\bf Z}^n \to Cl(\Sigma)$ has a left inverse $\phi :
Cl(\Sigma) \to {\bf Z}^n$.  Now, given $\alpha \in Cl(\Sigma)$,
consider the $S$-module $S(\alpha)$, which has the $Cl(\Sigma)$ grading
$S(\alpha)_\beta = S_{\alpha+\beta}$.  We can give $S(\alpha)$ a
grading by ${\bf Z}^n$ where $S(\alpha)_u = S_{u+\phi(\alpha)}$ for $u
\in {\bf Z}^n$ (and we are using the usual ${\bf Z}^n$ grading of
$S$).  By the above, this grading on $S(\alpha)$ gives a ${\bf
T}$-linearized sheaf ${\cal O}_{\bf P}(\alpha)$.

	If ${\cal F}$ is a quasi-coherent sheaf on ${\bf P}$, then \S3
of \cite{cox} implies that
\[ F = \bigoplus_{\alpha \in Cl(\Sigma)} H^0({\bf P}, {\cal F} \otimes
{\cal O}_{\bf P}(\alpha))\]
is a $Cl(\Sigma)$-graded $S$-module whose associated sheaf is ${\cal
F}$.  The module structure on $F$ comes from the natural isomorphism
$H^0({\bf P},{\cal O}_{\bf P}(\alpha)) \cong S_\alpha$.

	Now assume that ${\cal F}$ has a ${\bf T}$-linearization.
Then each ${\cal F} \otimes {\cal O}_{\bf P}(\alpha)$ has a natural
${\bf T}$-linearization, so that $H^0({\bf P},{\cal F} \otimes {\cal
O}_{\bf P}(\alpha))$ has a grading by $M$.  Then define a ${\bf Z}^n$
grading on $F$ by setting
\[ F_u = H^0({\bf P},{\cal F} \otimes {\cal
O}_{\bf P}(\alpha))_{u-\phi(\alpha)},\]
where $u$ maps to $\alpha$ under the map ${\bf Z}^n \to Cl(\Sigma)$.
Since $\phi$ is a homomorphism, it is easy to check that $F$ becomes a
${\bf Z}^n$-graded $S$-module and that this grading induces the given
${\bf T}$-linearization on ${\cal F}$.\hfill$\Box$

\begin{rem}
{\rm In Definition \ref{def.module}, we describe graded $S$ modules
which give the sheaves $\Omega_{\bf P}^p$.  The definition shows that
these modules have a canonical ${\bf Z}^n$ grading, which gives the
usual ${\bf T}$-linearization on $\Omega_{\bf P}^p$.  Note that we use
the ${\bf Z}^n$ grading in a crucial way in the proof of Proposition
\ref{prop.dminus1}. }
\end{rem}

\section{Local properties of differential forms}

Let $\Omega^p_{\bf P} $ denote the sheaf of $p$-differential forms of
Zariski on ${\bf P}$.  This means $\Omega^p_{\bf P} = j_*\Omega^p_W$,
where $W$ is the smooth part of ${\bf P}$, $\Omega^p_W$ is the usual
sheaf of $p$-forms on $W$, and $j : W \to {\bf P}$ is the natural
inclusion.  The sheaf $\Omega^p_{\bf P} $ has a canonical ${\bf
T}$-linearization which induces $M$-graded decompositions of sections
of $\Omega^p_{\bf P} $ over the ${\bf T}$-invariant the affine open
subsets ${\bf A}_\sigma \subset {\bf P}$.

\begin{opr}
{\rm Let $m$ be an element of $\check\sigma\cap M$, where $\sigma \in
\Sigma$. We denote by $\Gamma(m)$ the ${\bf C}$-subspace in $M_{\bf
C}$ generated by elements of the minimal face of $\check\sigma$
containing $m$. }
\end{opr}

\begin{prop}
\label{prop.A}
{\rm (\cite{dan1}, \S 4) } Let ${\bf A}_{\sigma} \subset {\bf P}$ be
the affine open corresponding to a $d$-dimensional cone $\sigma \in
\Sigma$.  Then the sections over ${\bf A}_\sigma$ of the {\bf
T}-linearized sheaf $\Omega_{\bf P}^p$ decompose into a direct sum of
$M$-homogeneous components as follows:
\[ \Omega_{{\bf A}_{\sigma}}^p : =
H^0({\bf A}_\sigma, \Omega_{{\bf P}}^p)
= \bigoplus_{m \in \check\sigma \cap M} \Lambda^p \Gamma(m). \]
\end{prop}

Besides $\Omega^p_{\bf P}$, we also have the sheaves $\Omega_{\bf
P}^p({\rm log}\, D)$ of differential
$p$-forms with logarithmic poles along $D = {\bf P} \setminus {\bf T}$
(see \S 15 of \cite{dan1}).  These sheaves
have the weight filtration
\[ {\cal W} : 0 \subset W_0  \Omega_{\bf P}^p({\rm log}\, D) \subset
W_1  \Omega_{\bf P}^p({\rm log}\, D) \subset \cdots \subset
W_p  \Omega_{\bf P}^p({\rm log}\, D) =  \Omega_{\bf P}^p({\rm log}\, D) \]
defined by
\[ W_k  \Omega_{\bf P}^p({\rm log}\, D) =
 \Omega_{\bf P}^{p-k} \wedge \Omega_{\bf P}^k({\rm log}\, D). \]
Note in particular that $W_0
\Omega_{\bf P}^p({\rm log}\, D) \cong \Omega_{\bf P}^p$ and $W_p
\Omega_{\bf P}^p({\rm log}\, D) \cong {\cal O}_{\bf P} \otimes
\Lambda^p M$.

	We have the following local description of the weight
filtration.

\begin{prop} {\rm (\cite{dan1}, \S 15.6)} If $\sigma$ is a
$d$-dimensional cone in $\Sigma$, then the sections over ${\bf
A}_\sigma$ of the {\bf T}-linearized sheaf $W_k \Omega_{\bf P}^p({\rm
log}\, D)$ decompose into a direct sum of $M$-homogeneous components
as follows:
\[ W_k \Omega_{{\bf A}_\sigma}^p({\rm log}\, D) := H^0({\bf A}_\sigma, W_k
\Omega_{\bf P}^p({\rm log}\, D)) = \bigoplus_{m \in \check\sigma
\cap M} \Lambda^{p-k}\Gamma(m) \wedge \Lambda^{k} M_{\bf C}. \]
\label{log.local}
\end{prop}

	The successive quotients of the weight filtration are described
using the Poincar\'e residue map.  Recall from Definition \ref{def.orbit}
that ${\bf P}_\tau$ is the Zariski closure of the ${\bf T}$-orbit of
${\bf P}$ corresponding to $\tau \in \Sigma$.

\begin{theo} {\rm (\cite{dan1}, \S 15.7)} For any integer
$k$ $(0 \leq k \leq p)$,
there is a short exact sequence
\[ 0 \rightarrow W_{k-1}  \Omega_{\bf P}^p({\rm log}\, D)
\rightarrow W_{k}  \Omega_{\bf P}^p({\rm log}\, D)
\stackrel{\rm Res}{\rightarrow}
\bigoplus_{\dim\tau = k}
\Omega_{{\bf P}_{\tau}}^p  \rightarrow 0 \]
where ``$\dim\tau = k$'' means the sum is over all $k$-dimensional cones
in $\tau \in \Sigma$, and ${\rm Res}$ is the Poincar\'e residue map.
\label{poin}
\end{theo}

This short exact sequence has a natural ${\bf T}$-action and splits
into $M$-homogeneous components.  Let's examine what happens over an
affine toric chart ${\bf A}_{\sigma}$, where $\sigma$ is a
$d$-dimensional cone in $\Sigma$.  Assume that the generators of
$\sigma$ are $\{e_1,\dots,e_d\}$.  Then, for any subset $\{ i_1,
\ldots, i_k \} \subset \{ 1, \dots, d \}$, we denote by ${\bf A}_{i_1
\dots i_k}$ the closed affine subvariety in ${\bf A}_{\sigma}$
corresponding to the cone $\tau_{i_1\dots i_k} = {\bf R}_{\geq 0}
e_{i_1} + \cdots + {\bf R}_{\geq 0} e_{i_k}$.  By \S 15.7 of
\cite{dan1}, we get the following local description of the residue
map:

\begin{prop}
\label{loc.poin}
Given $m \in \check\sigma\cap M$, let $\omega_m$ be an element of
$\Lambda^{p-k}\Gamma(m)$, and let $\omega'$ be an element of
$\Lambda^{k} M_{\bf C}$. Then the image of the $m$-homogeneous element
$\omega_m \wedge \omega' \in W_k \Omega_{{\bf A}_\sigma}^p({\rm log}\, D)$
under the residue map
\[ {\rm Res} : W_k \Omega_{{\bf A}_\sigma}^p({\rm log}\, D) \rightarrow
\bigoplus_{1 \leq i_1 < \cdots < i_k \leq d}
\Omega_{{\bf A}_{i_1 \dots i_k}}^{p-k} \]
is given by
\[ {\rm Res}(\omega_m \wedge \omega' )_{i_1 \dots i_k} =
\omega'(e_{i_1},
\ldots , e_{i_k}) \cdot \omega_m \in \Lambda^{p-k}\Gamma(m). \]
\end{prop}

\section{Globalization of the Poincar\'e residue map}

Let ${\cal L}$ be an ample {\bf T}-linearized invertible sheaf on the
complete toric variety ${\bf P}$, and let $\Delta = \Delta({\cal L})$
be as in Definition \ref{def.polytope}.  Tensoring by $\cal L$ the
short exact sequence in Theorem \ref{poin}, we obtain the exact
sequence
\[ 0 \rightarrow W_{k-1}  \Omega_{\bf P}^p({\rm log}\, D) \otimes {\cal L}
\rightarrow W_{k}  \Omega_{\bf P}^p({\rm log}\, D) \otimes {\cal L}
\stackrel{\rm Res}{\rightarrow}
\bigoplus_{\dim \tau = k}  \Omega_{{\bf P}_{\tau}}^{p-k}
\otimes {\cal L} \rightarrow 0. \]
The goal of this section is to give an explicit description of the
map of spaces of global sections
\begin{equation} \gamma :
H^0 ({\bf P}, W_{k}  \Omega_{\bf P}^p({\rm log}\, D) \otimes {\cal L})
\rightarrow \bigoplus_{\dim \tau = k}
H^0 ( {\bf P}, \Omega_{{\bf P}_{\tau}}^{p-k} \otimes {\cal L})
\label{global.poin}
\end{equation}
induced by the Poincar\'e residue map.

\begin{opr}
{\rm Let ${\cal L}$ be a ${\bf T}$-linearized ample invertible sheaf,
which determines the convex polytope in $\Delta = \Delta({\cal L})
\subset M_{\bf R}$.  For any $m \in M$, we denote by
$\Gamma_{\Delta}(m)$ the ${\bf C}$-subspace in $M_{\bf C}$ generated
by all vectors $s - s'$, where $s, s' \in \Delta_m$, and $\Delta_m$ is
the minimal face of $\Delta$ containing $m$. }
\label{glob.diff}
\end{opr}

First we notice the following properties:

\begin{prop}
\label{glob1}
The space of global sections
of ${\bf T}$-linearized sheaf $\Omega_{\bf P}^p({\rm log}\, D) \otimes
{\cal L}$
decomposes into a direct sum of $M$-homogeneous components as follows:
\[ H^0 ({\bf P}, W_{k}  \Omega_{\bf P}^p({\rm log}\, D) \otimes {\cal L}) =
\bigoplus_{m \in \Delta \cap M}
\Lambda^{p-k}\Gamma_{\Delta}(m) \wedge \Lambda^{k} M_{\bf C} \]
\end{prop}

\proof The statement follows from Proposition \ref{log.local} and
Proposition \ref{delta}. \hfill $\Box$

\begin{prop}
\label{glob2}
Let $\tau$ be a $k$-dimensional cone in $\Sigma$, $\Delta_{\tau}$ the
corresponding face of $\Delta$ of codimension $k$ $($see Proposition
\ref{ample}$)$.  Then the space of global sections of ${\bf
T}$-linearized sheaf $\Omega_{{\bf P}_{\tau}}^p \otimes {\cal L}$
decomposes into a direct sum of $M$-homogeneous components as follows:
\[ H^0 ({\bf P}, \Omega_{{\bf P}_{\tau}}^{p-k} \otimes {\cal L}) =
\bigoplus_{m \in \Delta_{\tau} \cap M}
\Lambda^{p-k}\Gamma_{\Delta_{\tau}}(m) \]
\end{prop}

\proof The statement follows from Proposition \ref{prop.A} and Proposition
\ref{delta}. \hfill $\Box$
\medskip

The linear mapping $\gamma$ of (\ref{global.poin}) is the direct sum
$\bigoplus_{\dim \tau = k} \gamma_{\tau}$, where
\[ \gamma_{\tau} \;:\;
H^0 ({\bf P}, W_{k}  \Omega_{\bf P}^p({\rm log}\, D) \otimes {\cal L})
\rightarrow
H^0 ( {\bf P}, \Omega_{{\bf P}_{\tau}}^{p-k} \otimes {\cal L}). \]
By Proposition \ref{loc.poin}, we can then describe $\gamma_\tau$ as follows.

\begin{prop}
\label{morph}
Let $\omega_m$ be an element of $\Lambda^{p-k}\Gamma_{\Delta}(m)$ and
$\omega'$ be an element of $\Lambda^{k} M_{\bf C}$.  Thus
$\omega_m \wedge \omega'$ is an  $m$-homogeneous element
in $H^0 ({\bf P}, W_{k}  \Omega_{\bf P}^p({\rm log}\, D) \otimes {\cal L})$.
Choose a $k$-dimensional cone $\tau$ with generators $e_{i_1},
\ldots, e_{i_k}$. Then
\begin{enumerate}
\item $\gamma_{\tau}(\omega_m \wedge \omega') = 0$ if
$m \not\in \Delta_{\tau}$;
\item $\gamma_{\tau}(\omega_m \wedge \omega') =
\omega'(e_{i_1}, \ldots , e_{i_k}) \cdot \omega_m$
if $m \in \Delta_{\tau}$.
\end{enumerate}
\end{prop}

\section{A generalized theorem of Bott-Steenbrink-Danilov}

In Theorem 7.5.2 of \cite{dan1}, Danilov formulated without proof the
following vanishing theorem generalizing for complete toric varieties
the well-known theorem of Bott and Steenbrink:

\begin{theo}
Let $\cal L$ be an ample invertible sheaf on $\bf P$. Then
for any $p \geq 0$ and $ i > 0$, one has
\[ H^i({\bf P}, \Omega^p_{{\bf P}} \otimes {\cal L}) =0. \]
\end{theo}

We prove now for simplicial toric varieties a more general vanishing
theorem:

\begin{theo}
\label{theo.bott}
Let $\cal L$ be an ample invertible sheaf on a complete simplicial
toric variety $\bf P$. Then for any $p \geq 0$, $k \geq 0$ and $ i >
0$, one has
\[ H^i({\bf P}, W_k\Omega^p_{{\bf P}}({\rm log}\, D) \otimes {\cal L})
= 0. \]
\end{theo}

\proof We prove this theorem using induction on $p-k$. For $p-k =0$,
the sheaf $W_p\Omega^p_{{\bf P}}({\rm log}\, D) \otimes {\cal L} =
\Omega_{\bf P}^p({\rm log}\, D)\otimes{\cal L} = \Lambda^p M_{\bf
C}\otimes {\cal L}$ is the direct sum of ${ d \choose p} $ copies of
$\cal L$. Thus, the vanishing property for $W_p\Omega^p_{{\bf P}}({\rm
log}\, D) \otimes {\cal L}$ is implied by the following general
vanishing property for the ample invertible sheaf $\cal L$.

\begin{prop} {\rm (\cite{dan1}, \S 7.3)}
Let ${\cal L}$ be an ample invertible sheaf on a complete toric variety
${\bf P}$. Then
\[ H^i({\bf P}, {\cal L}) = 0\;\; {\rm for}\;\;i > 0. \]
\end{prop}

On the other hand, for any $k$ ($0 \leq k \leq p$) we can apply the
induction assumption to $W_{k} \Omega_{\bf P}^p({\rm log}\, D) \otimes
{\cal L}$ and $\Omega_{{\bf P}_{\tau}}^p \otimes {\cal L}$ appearing
in the short exact sequence
\[ 0 \rightarrow W_{k-1}  \Omega_{\bf P}^p({\rm log}\, D) \otimes {\cal L}
\rightarrow W_{k}  \Omega_{\bf P}^p({\rm log}\, D) \otimes {\cal L}
\stackrel{\rm Res}{\rightarrow}
\bigoplus_{\dim\tau = k }  \Omega_{{\bf P}_{\tau}}^{p-k} \otimes {\cal L}
\rightarrow 0. \]
The required vanishing properties of  $W_{k-1}
\Omega_{\bf P}^p({\rm log}\, D) \otimes {\cal L}$ follows now from   the
following lemma.

\begin{lem}
The mapping
\[ \gamma : H^0 ({\bf P}, W_{k}  \Omega_{\bf P}^p({\rm log}\, D)
\otimes {\cal L})  \rightarrow
\bigoplus_{\dim\tau = k}
H^0 ( {\bf P}, \Omega_{{\bf P}_{\tau}}^{p-k} \otimes {\cal L}) \]
is surjective.
\end{lem}

\proof Since ${\cal L}$ is ample, there exists a one-to-one
correspondence between $i$-di\-men\-sional cones of $\Sigma$ and
$(d-i)$-dimensional faces of the convex polytope $\Delta =
\Delta({\cal L})$ (see Proposition \ref{ample}).  Choose a
$k$-dimensional cone $\tau_0 \in \Sigma$.  We know from Proposition
\ref{glob2} that the $m$-homogeneous component of $H^0 ( {\bf
P}_{\tau_0}, \Omega_{{\bf P}_{\tau_0}}^{p-k} \otimes {\cal L})$ is
non-zero only if $m \in \Delta_{\tau_0}$, and when the latter
holds, Definition \ref{glob.diff} shows that the $m$-component is
determined by the minimal face $\Delta_{\tau_1}$ of $\Delta_{\tau_0}$
containing $m$.  This means that $\tau_1 \in \Sigma$ is a cone of
dimension $c \geq k$ containing $\tau_0$.

	Fix such a lattice point $m$ and the corresponding $\tau_0
\subset \tau_1$.  We can assume that $\tau_1$ is a face of a $d$-dimensional
cone $\sigma \in \Sigma$ with generators $e_1, \ldots, e_d \in N$
such that $e_1, \ldots, e_c$ are generators of $\tau_1$ and $e_1,
\ldots, e_k$ are generators of $\tau_0$.  To describe
$\Gamma_{\Delta_{\tau_0}}(m)$, suppose that $\Delta$ is defined by
inequalities $\langle m',e_i\rangle \ge -a_i$ for $1 \le i \le n$.
Then $\Delta_{\tau_1} \subset \Delta$ is the face obtained by
requiring $\langle m',e_i\rangle = -a_i$ for $1 \le i \le c$, which
implies that the subspace $\Gamma_{\Delta_{\tau_0}}(m) \subset M_{\bf
C}$ is defined $\langle m',e_i\rangle = 0$ for $1 \le i \le c$.

	Now let $h_1, \ldots, h_d \in M_{\bf Q}$ form the dual basis
to the basis $e_1, \ldots, e_d$ of $N_{\bf Q}$.  It follows
immediately that $h_{c+1},\dots,h_d$ form a basis of
$\Gamma_{\Delta_{\tau_0}}(m)$.  Thus, by Proposition
\ref{glob2}, the $m$-homogeneous component of $H^0 ( {\bf P},
\Omega_{{\bf P}_{\tau_0}}^{p-k} \otimes {\cal L})$ has a basis
consisting of $(p-k)$-vectors
\[ \omega_{j_1 \ldots j_{p-k}} = h_{j_1} \wedge \cdots \wedge h_{j_{p-k}} \]
where $\{ j_1, \ldots, j_{p-k} \}$ is a subset of
$\{ c+1, \ldots, d\}$. For any such a $(p-k)$-vector
$\omega_{j_1 \ldots j_{p-k}}$, the $p$-vector
\[ h_{1} \wedge \cdots
\wedge h_{{k}} \wedge \omega_{j_1 \ldots j_{p-k}} \]
defines an element in the $m$-homogeneous component of
the space
$H^0 ({\bf P}, W_{k}  \Omega_{\bf P}^p({\rm log}\, D) \otimes {\cal L})$
by Proposition \ref{glob1}.   Furthermore, Proposition \ref{morph}
shows that
\[ \gamma_{\tau_0}(h_{1} \wedge \cdots \wedge h_{{k}} \wedge
\omega_{j_1 \ldots j_{p-k}}) = h_1\wedge\cdots \wedge h_k (e_1\wedge
\cdots \wedge e_k)\cdot\omega_{j_1 \ldots j_{p-k}} = \omega_{j_1
\ldots j_{p-k}}.\]
Thus, to prove the lemma, it suffices to show that for $k$-dimensional
cones $\tau \in \Sigma$, we have
\[ \gamma_{\tau}(h_{1} \wedge \cdots \wedge h_{k} \wedge
\omega_{j_1 \ldots j_{p-k}}) = 0 \; \; {\rm for}\; \tau \neq \tau_0.  \]

However, by Proposition \ref{morph},
$\gamma_{\tau}(h_{1} \wedge \cdots \wedge h_{{k}} \wedge
\omega_{j_1 \ldots j_{p-k}}) = 0$ for $k$-dimensional cones $\tau \in
\Sigma$ such that $m \notin \Delta_{\tau}$.  Since $\Delta_{\tau_1}$
is the minimal face containing $m$, the condition $m \notin
\Delta_\tau$ holds whenever $\tau$ is not a face of $\tau_1$.  It
remains to see what happens when $\tau$ is a $k$-dimensional face
of $\tau_1$.  In this case, the generators of $\tau$ are $e_{i_1},
\ldots , e_{i_k}$, where $\{ i_1, \ldots, i_k \} \subset \{ 1, \ldots,
c\}$.  Since the $e_i$ are dual to the $h_i$, the value $h_{1} \wedge
\cdots \wedge h_{k}(e_{i_1}, \ldots , e_{i_k})$ is nonzero only if
$\tau = \tau_0$.  From Proposition \ref{morph}, it follows that
$\gamma_{\tau}(h_{1} \wedge \cdots \wedge h_{{k}} \wedge \omega_{j_1
\ldots j_{p-k}}) = 0$ when $\tau \ne \tau_0$. \hfill $\Box$

\section{Differential forms and graded $S$-modules}

As usual, $\Omega_{\bf P}^p$ denotes the sheaf of Zariski
differential $p$-forms on a complete simplicial toric variety ${\bf P}
= {\bf P}_\Sigma$.  We will study these sheaves using certain graded
modules over the polynomial ring $S = {\bf C}[z_1,\dots,z_n]$, which
is graded by $Cl(\Sigma)$.  Given a Cartier divisor $X \subset {\bf
P}$, our goal is to describe $H^0({\bf P},\Omega^p_{\bf P}(X))$ in
terms of $S$.

	Recall that the fan $\Sigma$ lies in $N_{\bf R} \cong {\bf
R}^d$ and that $M = {\rm Hom}(N,{\bf Z})$.  Also recall that
$e_1,\dots,e_n$ are the generators of the 1-dimensional cones of
$\Sigma$.

\begin{opr} {\rm Given $p$ between $0$ and $d$, define
$\widehat{\Omega}^p_S$ by the exact sequence of graded $S$-modules:
\[ 0 \to \widehat{\Omega}_S^p \to S \otimes \Lambda^p M
\stackrel{\gamma}\rightarrow \bigoplus_{i=1}^n
(S/z_i S)\otimes\Lambda^{p-1} (M\cap e_i^\perp), \]
where the $i^{\rm th}$ component of $\gamma$ is $\gamma_i(g\otimes
\omega) = g \bmod z_i \otimes \langle e_i,\omega\rangle$.  A careful
description of the interior product $\langle e_i,w\rangle \in
\Lambda^{p-1}(M\cap e_i^\perp)$ may be found in \S3.2 of \cite{oda}.}
\label{def.module}
\end{opr}

	By \cite{cox}, we know that every finitely generated
graded $S$-module $F$ gives a coherent sheaf $\widetilde F$ on
${\bf P}$.

\begin{lem} The sheaf $\widetilde{\widehat{\Omega}}{}^p_S$ on ${\bf P}$
associated to the graded $S$-module $\widehat{\Omega}^p_S$ is naturally
isomorphic to $\Omega^p_{\bf P}$.
\label{lem.tilde}
\end{lem}

\proof From \S3 of \cite{cox}, we get $\widetilde{S} = {\cal
O}_{\bf P}$.  The primitive element $e_i$ generates the cone $\rho_i
\in \Sigma$, which by Definition \ref{def.orbit} gives the ${\bf
T}$-invariant divisor $D_i = {\bf P}_{\rho_i} \subset {\bf P}$.  By
the example before Corollary 3.8 of \cite{cox}, the ideal $z_i S$
gives the ideal sheaf of $D_i$.  Since $F \mapsto
\widetilde{F}$ is exact (see Proposition 3.1 of \cite{cox}), it
follows that $S/z_i S$ gives the sheaf ${\cal O}_{D_i}$.  Then the
exactness of $F \mapsto \widetilde{F}$, applied to the sequence
defining $\Omega_S^p$, gives the exact sequence of sheaves:
\[ 0 \to \widetilde{\widehat{\Omega}}{}^p_S \to {\cal O}_{\bf P}\otimes
\Lambda^p M \stackrel{\gamma}\to \bigoplus_{i=1}^n {\cal O}_{D_i}\otimes
\Lambda^{p-1} (M\cap e_i^\perp)\ .\]
By Theorem 3.6 of \cite{oda}, we can identify
$\widetilde{\widehat{\Omega}}{}^p_S$ with $\Omega^p_{\bf P}$.
\hfill$\Box$
\medskip

	We next discuss the twists of the $\Omega^p_{\bf P}$.  First
recall that if $\beta \in Cl(\Sigma)$ and $F$ is a graded
$S$-module, then $F(\beta)$ is the graded $S$-module defined by
$F(\beta)_\gamma = F_{\beta+\gamma}$ for $\gamma \in Cl(\Sigma)$.

\begin{opr} {\rm Given $\beta \in Cl(\Sigma)$, the sheaf on ${\bf P}$
associated to the graded $S$-module $\widehat{\Omega}_S^p(\beta)$ is denoted
$\Omega_{\bf P}^p(\beta)$.}
\end{opr}

\begin{rem} {\rm If ${\cal O}_{\bf P}(\beta)$ is the sheaf associated
to $S(\beta)$, then there is a natural isomorphism $\Omega^p_{\bf
P}(\beta) \cong \Omega^p_{\bf P}\otimes {\cal O}_{\bf P}(\beta)$
whenever $\beta$ is the class of a Cartier divisor.  However, when
$\beta$ is not Cartier, the sheaves $\Omega^p_{\bf P}(\beta)$ and
$\Omega^p_{\bf P}\otimes {\cal O}_{\bf P}(\beta)$ may be nonisomorphic.}
\end{rem}

\begin{prop}
\label{prop.sections}
For any divisor class $\beta \in Cl(\Sigma)$, there is a natural
isomorphism
\[ H^0({\bf P},\Omega^p_{\bf P}(\beta)) \cong (\widehat{\Omega}^p_S)_\beta.\]
\end{prop}

\proof The sequence in Definition \ref{def.module} remains exact after
shifting by $\beta$.  Then, taking the associated sheaves on ${\bf P}$,
we get an exact sequence
\[ 0 \to \Omega^p_{\bf P}(\beta) \to {\cal O}_{\bf P}(\beta)\otimes
\Lambda^p M \stackrel{\gamma}\to \bigoplus_{i=1}^n {\cal
O}_{D_i}(\beta) \otimes \Lambda^{p-1} (M\cap e_i^\perp)\ .\]
Since taking global sections is left exact, we get the exact sequence
\[0 \to H^0({\bf P},\Omega^p_S(\beta)) \to H^0({\bf P},{\cal O}_{\bf
P}(\beta)) \otimes
\Lambda^p M \stackrel{\gamma}\to \bigoplus_{i=1}^n H^0({\bf P},{\cal
O}_{D_i}(\beta)) \otimes \Lambda^{p-1} (M\cap e_i^\perp)\ .\]
Using the natural isomorphism $H^0({\bf P},{\cal O}_{\bf P}(\beta))
\cong S_\beta$ from Proposition 1.1 of \cite{cox}, we get
\begin{equation} 0 \to H^0({\bf P},\Omega^p_S(\beta)) \to S_\beta \otimes
\Lambda^p M \stackrel{\gamma}\to \bigoplus_{i=1}^n H^0({\bf P},{\cal
O}_{D_i}(\beta)) \otimes \Lambda^{p-1} (M\cap e_i^\perp).
\label{eq.first}
\end{equation}
However, since $z_i$ vanishes on $D_i$, the map $\gamma$ factors
\[S_\beta \otimes \Lambda^p M \to \bigoplus_{i=1}^n (S/z_i S)_\beta
\otimes \Lambda^{p-1} (M\cap e_i^\perp) \to
\bigoplus_{i=1}^n H^0({\bf P},{\cal
O}_{D_i}(\beta)) \otimes \Lambda^{p-1} (M\cap e_i^\perp).\]
Assume for the moment that $(S/z_i S)_\beta \to
H^0({\bf P},{\cal O}_{D_i}(\beta))$ is injective.  Then
(\ref{eq.first}) gives the exact sequence
\[ 0 \to H^0({\bf P},\Omega^p_S(\beta)) \to S_\beta \otimes
\Lambda^p M \stackrel{\gamma}\to \bigoplus_{i=1}^n (S/z_i S)_\beta
\otimes \Lambda^{p-1} (M\cap e_i^\perp),\]
and the isomorphism $H^0({\bf P},\Omega^p_S(\beta)) \cong
(\widehat{\Omega}^p_S)_\beta$ follows immediately from Definition
\ref{def.module}.

	It remains to show that the natural map $(S/z_i S)_\beta \to
H^0({\bf P},{\cal O}_{D_i}(\beta))$ is injective.  Let $\sigma$ be a
cone of $\Sigma$ containing $e_i$, and let ${\bf A}_\sigma \subset
{\bf P}$ be the corresponding affine toric variety.  Then it suffices
to show that the $(S/z_i S)_\beta \to H^0({\bf A}_\sigma,{\cal
O}_{D_i}(\beta))$ is injective.  However, as explained in \S3 of
\cite{cox}, we have
\[ H^0({\bf A}_\sigma,{\cal O}_{D_i}(\beta)) \cong ((S/z_i S)(\beta)
\otimes S_\sigma)_0 = ((S/z_iS)\otimes S_\sigma)_\beta,\]
where $S_\sigma$ is the localization of $S$ at $\widehat{z}_\sigma =
\prod_{e_j \notin \sigma} z_j$.  Thus we need to show that
\[ (S/z_i S)_\beta \to ((S/z_i S)\otimes S_\sigma)_\beta\]
is injective.  This follows easily since $z_i$ doesn't divide
$\widehat{z}_\sigma$, and the proposition is proved. \hfill $\Box$
\medskip

	We next want to relate $\widehat{\Omega}^p_S$ to the usual
module $\Omega^p_S$ of $p$-forms in $dz_1,\dots,dz_n$.  First note
that $\Omega^\cdot_S$ can be given the structure of a graded
$S$-algebra by declaring that $dz_i$ and $z_i$ have the same degree in
$Cl(\Sigma)$.  This means that if $z_i \in S_{\beta_i}$, then we have
an isomorphism of graded $S$-modules
\[ \Omega^p_S \cong \bigoplus_{1 \le i_1 < \cdots < i_p \le n}
S(-\beta_{i_1} - \cdots - \beta_{i_p}).\]

\begin{lem} There is a natural inclusion $\widehat{\Omega}^p_S \subset
\Omega^p_S$ of graded $S$-modules.
\end{lem}

\proof  Consider the diagram
\begin{equation}
\begin{array}{ccccccc}
0 & \to & \widehat{\Omega}^p_S & \to & S\otimes\Lambda^p M &
\stackrel{\gamma}\to &
\bigoplus_{i=1}^n (S/z_i S)\otimes\Lambda^{p-1} M \\
& & & & \downarrow & & \downarrow \\
0 & \to & \Omega^p_S & \to & S\otimes\Lambda^p {\bf Z}^n &
\stackrel{\delta}\to &
\bigoplus_{i=1}^n (S/z_i S)\otimes\Lambda^{p-1} {\bf Z}^n
\label{dia.omega}
\end{array}
\end{equation}
The first row is exact by the definition of $\widehat{\Omega}^p_S$ and
the inclusion $M\cap e_i^\perp \subset M$.  To understand the second row,
let $h_1,\dots,h_n$ be the standard basis of ${\bf Z}^n$.  Then the map
$\Omega^p_S \to S\otimes\Lambda^p{\bf Z}^n$ is defined by
\[ dz_{i_1}\wedge\cdots\wedge dz_{i_p} \mapsto z_{i_1}\cdots
z_{i_d}\otimes h_{i_1}\wedge\cdots\wedge h_{i_d},\]
and for $g\otimes\omega \in S\otimes\Lambda^p{\bf Z}^n$, the $i^{\rm
th}$ component of $\delta(g\otimes\omega)$ is $\delta_i(g\otimes\omega)
= g \bmod z_i \otimes \langle h_i^*,\omega\rangle$, where
$h_1^*,\dots,h_n^*$ is the dual basis to $h_1,\dots,h_n$.  It is easy
to check that the sequence on the bottom is exact.

	The vertical maps in (\ref{dia.omega}) are induced by the map
$\alpha : M \to {\bf Z}^n$ from Definition \ref{diag}, and since the
dual $\alpha^* : {\bf Z}^n \to N$ maps $h_i^*$ to $e_i$, it follows
that diagram (\ref{dia.omega}) commutes.  This gives the desired
inclusion. \hfill $\Box$
\medskip

	The final step is to relate $\widehat{\Omega}^p_S$ to rational
$p$-forms on ${\bf P}$ with poles on a hypersurface $X \subset {\bf
P}$.  We will assume that $X$ is defined by $f = 0$ for some $f \in
S_\beta$.  Thus $\beta \in Cl(\Sigma)$ is the class of $X$, and we
will also assume that $X$ is a Cartier divisor.  Then $\Omega^p_{\bf
P}(X) = \Omega^p_{\bf P} \otimes {\cal O}_{\bf P}(X)$.  Recall that
local sections of $\Omega^p_{\bf P}(X)$ are rational $p$-forms which
become holomorphic when multiplied by the local equation of $X$.

\begin{prop}
\label{prop.pform}
If $X$ is a Cartier divisor on ${\bf P}$ defined by $f =
0$ for $f \in S_\beta$, then
\[ H^0({\bf P},\Omega^p_{\bf P}(X)) =
\Bigl\{ {\omega \over f} : \omega \in
(\widehat{\Omega}^p_S)_\beta \Bigr\}.\]
\end{prop}

\proof We first observe that multiplication by $1/f$ gives an
isomorphism ${\cal O}_{\bf P}(\beta) \cong {\cal O}_{\bf P}(X)$.  To
see this, we will work locally on an affine open ${\bf A}_\sigma
\subset {\bf P}$, where $\sigma \in \Sigma$.  Recall that ${\bf
A}_\sigma = {\rm Spec}\,(S_\sigma)_0$, where $S_\sigma$ is the
localization of $S$ at $\widehat{z}_\sigma = \prod_{e_i \notin \sigma}
z_i$.  Since $X$ is Cartier, Lemma 3.4 of \cite{cox} shows that there
is a monomial $z^D \in S_\beta$ which is invertible in $S_\sigma$.  It
follows that $f/z^D = 0$ is the local equation of $X$ on ${\bf
A}_\sigma$.  Hence
\[ H^0({\bf A}_\sigma,{\cal O}_{\bf P}(X)) = {1 \over f/z^D} \cdot
H^0({\bf A}_\sigma,{\cal O}_{\bf P}) = {1 \over f/z^D} \cdot
(S_\sigma)_0 = {1 \over f} \cdot (S_\sigma)_\beta ,\]
so that $(1/f)S$ is the graded $S$-module that gives ${\cal
O}_{\bf P}(X)$.

	Thus $1/f : S(\beta) \to (1/f)S$ is a graded isomorphism which
induces ${\cal O}_{\bf P}(\beta) \cong {\cal O}_{\bf P}(X)$.  This
proves that we have an isomorphism $\Omega^p_{\bf P}(\beta) \cong
\Omega^p_{\bf P}(X)$ which is given by multiplication by $1/f$
when we represent each sheaf by a graded $S$-module.  Then the desired
result follows immediately from Proposition \ref{prop.sections}.
\hfill $\Box$

\begin{rem} {\rm There is a direct way of seeing that $\omega/f$
descends to $p$-form on ${\bf P}$.  In diagram (\ref{dia.omega}), we
can think of $S\otimes\Lambda^p{\bf Z}^n$ as $\Omega^p_S({\rm
log}\,\tilde{D})$, where $\tilde{D}$ is the union of the coordinate
hyperplanes.  With this interpretation, the standard basis of ${\bf
Z}^n$ is $dz_i/z_i$, and the map $\Omega^p_S \to S\otimes
\Lambda^p{\bf Z}^n$ is given by $dz_i \mapsto z_i\otimes dz_i/z_i$.
Now fix a basis $m_1,\dots,m_d$ of $M$ and let
\[ t_j = \prod_{i=1}^n z_i^{\langle m_j,e_i\rangle}\]
for $1 \le j \le d$.  The $t_j$ are invariant under the group ${\bf D}
= {\bf D}(\Sigma)$ and hence descend to rational functions on ${\bf
P}$ (in fact, they are coordinates for the torus ${\bf T}
\subset {\bf P}$).  Note that $dt_j/t_j = \sum_{i=1}^n\langle
m_j,e_i\rangle dz_i/z_i = \alpha(m_j)$, where $\alpha : M \to {\bf
Z}^n$ is the map from Definition \ref{diag}.  Then, in
$S\otimes \Lambda^p{\bf Z}^n$, we can write
\begin{eqnarray*} \omega & = & \sum_{j_1 < \cdots < j_d}
g_{{j_1}\cdots{j_d}} \otimes \alpha(m_{j_1})\wedge\cdots\wedge
\alpha(m_{j_d})\\
& = & \sum_{j_1 < \cdots < j_d}
g_{{j_1}\cdots{j_d}} \otimes dt_{j_1}/t_{j_1}\wedge\cdots\wedge
dt_{j_d}/t_{j_d}
\end{eqnarray*}
where $g_{{j_1}\cdots{j_d}} \in S_\beta$.  Thus $\omega/f$ can be
written
\[{\omega \over f} = \sum_{j_1 < \cdots < j_d}
{g_{{j_1}\cdots{j_d}} \over f} \otimes dt_{j_1}/t_{j_1}\wedge\cdots\wedge
dt_{j_d}/t_{j_d}.\]
Since $g_{{j_1}\cdots{j_d}}$ and $f$ have the same degree,
$\omega/f$ descends to a rational $p$-form on ${\bf P}$.}
\label{rem.dz}
\end{rem}

\section{Differential forms of degree $d$ and $d-1$}

	In this section, we will find module generators for
$\widehat{\Omega}^d_S$ and $\widehat{\Omega}^{d-1}_S$, where $d$ is
the dimension of the complete simplicial toric variety ${\bf P}$.
This will enable us to give an explicit description of $H^0({\bf
P},\Omega^p_{\bf P}(X))$ for $p = d$ and $d-1$.

\begin{opr} {\rm Fix an integer basis $m_1,\dots,m_d$ for the lattice
$M$.  Then, given a subset $I = \{i_i,\dots,i_d\} \subset \{1,\dots,n\}$
consisting of $d$ elements, define
\[ \det(e_I) = \det(\langle
m_j,e_{i_k}\rangle_{\scriptscriptstyle{1 \le j,k \le d}}). \]
We also define $dz_I = dz_{i_1}\wedge\cdots\wedge dz_{i_d}$ and
$\widehat{z}_I = \textstyle{\prod_{i \notin I}} z_i$.}
\label{def.eI}
\end{opr}

\begin{rem} {\rm Although $\det(e_I)$ and $dz_I$ depend on how the
elements of $I$ are ordered, their product $\det(e_I) dz_I$ does not.}
\end{rem}

\begin{opr}
\label{def.omegao}
{\rm We define the $d$-form $\Omega_0 \in \Omega^d_S$ by the
formula
\[\Omega_0 = \sum_{|I| = d} \det(e_I) \widehat{z}_I dz_I,\]
where the sum is over all $d$ element subsets $I \subset
\{1,\dots,n\}$.}
\end{opr}

\begin{exam} {\rm Suppose that ${\bf P}$ is the weighted projective
space ${\bf P}(w_1,\dots,w_{d+1})$.  Then one can check that up to a
nonzero constant (which depends on the basis of $M$ chosen in
Definition \ref{def.eI}), we have
\[ \Omega_0 = \sum_{i=1}^{d+1} (-1)^i w_i z_i dz_1\wedge\cdots
\wedge \widehat{dz_i} \wedge\cdots \wedge dz_{d+1}.\]
This form appears in 2.1.3 of \cite{dolgach}, where it is denoted
$\Delta(dT_0\wedge\cdots\wedge dT_r)$.  Also, when ${\bf P} = {\bf
P}^d$, we have $w_i = 1$ for all $i$, and we recover the form $\Omega$
in Corollary 2.11 of \cite{griff1}.}
\end{exam}

\begin{prop} $\widehat{\Omega}^d_S \subset \Omega^d_S$ is a free
$S$-module of rank 1 generated by $\Omega_0 \in \Omega^d_S$.
\label{prop.dmodule}
\end{prop}

\proof Using the basis $m_1,\dots,m_d$ of $M$, we see that
$S\otimes\Lambda^d M$ is free of rank 1 with generator
$m_1\wedge\cdots\wedge m_d$.  Furthermore, since $\langle
e_i,m_1\wedge\cdots\wedge m_d\rangle$ is nonzero in the rank 1 {\bf
Z}-module $\Lambda^{d-1}(M\cap e_i^\perp)$, the definition of
$\widehat{\Omega}^d_S$ shows that for $A \in S$, $A\,
m_1\wedge\cdots\wedge m_d$ lies in $\widehat{\Omega}^d_S$ if and only
if $A$ is divisible by $z_1\cdots z_n$.  Thus $\widehat{\Omega}^d_S
\subset S\otimes\Lambda^d M$ is free of rank 1 with $z_1\dots z_n\,
m_1\wedge\cdots\wedge m_d$ as generator.

	As in Remark \ref{rem.dz}, we will denote the standard basis
of ${\bf Z}^n$ by $dz_1/z_1,\dots,dz_n/z_n$.  Then, using the map
$\alpha : M \to {\bf Z}^n$, we have $\alpha(m) = \sum_{i=1}^n \langle
m,e_i\rangle dz_i/z_i$.  Thus, inside $S\otimes\Lambda^d{\bf Z}^n$,
the generator $z_1\cdots z_n\,m_1\wedge\cdots\wedge m_d$ of
$\widehat{\Omega}^d_S$ is
\begin{eqnarray*}
z_1\cdots z_n \textstyle{\left(\sum_{i=1}^n \langle
m_1,e_i\rangle {dz_i\over z_i}\right) \wedge \cdots \wedge \left(\sum_{i=1}^n
\langle m_d,e_i\rangle {dz_i\over z_i}\right)} &  = &
z_1\cdots z_n \sum_{|I| = d}
\det(e_I) (\textstyle{\prod_{i \in I} z_i^{-1}}) dz_I\\
	& = & \sum_{|I| = d} \det(e_I) \widehat{z}_I dz_I\ =\ \Omega_0,
\end{eqnarray*}
and the proposition is proved. \hfill$\Box$

\begin{rem}
\label{rem.omegad}
{\rm If $z_i \in S_{\beta_i}$ for $1 \le i \le n$, then we
set $\beta_0 = \sum_{i=1}^n \beta_i.$ Since $\Omega_0$ has degree
$\beta_0$, the above proposition gives an isomorphism of graded
$S$-modules
\[ \widehat{\Omega}^d_S \cong S(-\beta_0).\]
Furthermore, since $\beta_0 \in Cl(\Sigma)$ is the class of ${\bf P}
\backslash {\bf T} = \sum_{i=1}^n D_i$, we get the well-known
isomorphism of sheaves
\[ \Omega^d_{\bf P} \cong {\cal O}_{\bf
P}\bigl(-\textstyle{\sum_{i=1}^n} D_i\bigr).\]}
\end{rem}

	If we combine Propositions \ref{prop.dmodule} and
\ref{prop.pform}, we get the following way of representing $d$-forms
on $P$ with poles on $X$.

\begin{theo}
\label{theo.dform}
Let $X \subset {\bf P}$ be a Cartier divisor defined by
$f = 0$, where $f \in S_\beta$.  Then
\[ H^0({\bf P},\Omega^d_{\bf P}(X)) =
\Bigl\{ {A\Omega_0 \over f} : A \in S_{\beta-\beta_0}\Bigr\},\]
where $\beta_0 = \sum_{i=1}^n \beta_i$ and $z_i \in S_{\beta_i}$.
\end{theo}

	We next describe generators for $\widehat{\Omega}^{d-1}_S$.

\begin{opr}
\label{def.omegai}
{\rm Given $i$ between $1$ and $n$, we define the $(d-1)$-form
$\Omega_i \in \Omega^{d-1}_S$ by the formula
\[ \Omega_i = \sum_{\scriptstyle{ |J| = d-1 \atop i \notin J}}
\det(e_{J\cup\{i\}}) \widehat{z}_{J\cup\{i\}} dz_J, \]
where $\det(e_{J\cup\{i\}})$ is computed by ordering the elements of
$J\cup\{i\}$ so that $i$ is {\em first\/} (this ensures that
$\det(e_{J\cup\{i\}}) dz_J$ is well-defined).}
\end{opr}

\begin{exam} {\rm When ${\bf P} = {\bf P}^1\times{\bf P}^1$, we have
$e_1 = (1,0)$, $e_2 = -e_1$, $e_3 = (1,0)$ and $e_4 = -e_3$, and one
can check that
\label{exam.p1p1}
\begin{eqnarray*}
\Omega_1 & = & z_2(z_4dz_3 - z_3dz_4) \\
\Omega_2 & = & -z_1(z_4dz_3 - z_3dz_4) \\
\Omega_3 & = & z_3(z_1dz_2 - z_2dz_1) \\
\Omega_4 & = & -z_4(z_1dz_2 - z_2dz_1)
\end{eqnarray*}
In this case, one can show that $\widehat{\Omega}^1_S$ is the
submodule of $\Omega^1_S$ generated by $z_1dz_2-z_2dz_1$ and $z_4dz_3
- z_3dz_4$ (this will also follow from Proposition \ref{prop.dminus1}
below).}
\end{exam}

	As the above example indicates, the $\Omega_i$ may be
multiples of other forms.  This happens whenever there is a $j$ such
that $e_j = -e_i$.  In this situation, note that $\det(e_{J\cup\{i\}})
= 0$ when $j \in J$ since the matrix will have one column which is the
negative of another.  But when $j \notin J$, then $z_j$ divides
$\widehat{z}_{J\cup\{i\}}$.  It follows that $e_j = -e_i$ implies
\[ \Omega_i = z_j \sum_{\scriptstyle{|J| = d-1 \atop i,j \notin J}}
\det(e_{J\cup\{i\}}) \widehat{z}_{J\cup\{i,j\}} dz_J. \]
Hence we get the following definition.

\begin{opr} {\rm Let
\begin{eqnarray*}
	{\cal I}_0 & = & \bigl\{ i : -e_i \notin
	\{e_1,\dots,e_n\}\bigr\} \\
	{\cal I}_1 & = & \bigl\{ (i,j) : i < j,\ e_i = -e_j\bigr\}
\end{eqnarray*}
Furthermore, for $(i,j) \in {\cal I}_1$, we define the $(d-1)$-form
\[ \Omega_{ij} = \sum_{\scriptstyle{|J| = d-1 \atop i,j \notin J}}
\det(e_{J\cup\{i\}}) \widehat{z}_{J\cup\{i,j\}} dz_J. \]}
\end{opr}

\begin{rem} {\rm For $(i,j) \in {\cal I}_1$, we have $\Omega_i =
z_j\Omega_{ij}$, and since $\det(e_{J\cup\{i\}}) =
-\det(e_{J\cup\{j\}})$ (remember $e_j = -e_i$), we get $\Omega_j = -
z_i\Omega_{ij}$ as in Example \ref{exam.p1p1}.}
\end{rem}

\begin{prop} With the above notation, $\widehat{\Omega}^{d-1}_S$ is
the submodule of $\Omega^{d-1}_S$ generated by the $\Omega_i$ for $i \in
{\cal I}_0$ and the $\Omega_{ij}$ for $(i,j) \in {\cal I}_1$.
\label{prop.dminus1}
\end{prop}

\proof We begin with the exact sequence
\begin{equation}
 0 \to \widehat{\Omega}^{d-1}_S \to S\otimes\Lambda^{d-1}
\stackrel{\gamma}\to \bigoplus_{i=1}^n (S/z_i
S)\otimes\Lambda^{d-2}M
\label{seq.dminus1}
\end{equation}
which follows from the definition of $\widehat{\Omega}^{d-1}_S$ and
the inclusion $M\cap e_i^\perp \subset M$.  If $m_1,\dots,m_d$ is an
integer basis of $M$, then we get the bases
\[ \begin{array}{rl}
\hbox{basis of $\Lambda^{d-1}M$}: & \omega_j = (-1)^{j+1}m_1\wedge
\cdots \wedge \widehat{m_j} \wedge \cdots \wedge m_d,\quad 0 \le j \le
d \\
\hbox{basis of $\Lambda^{d-2}M$}: & \omega_{jk} = (-1)^{j+k}m_1\wedge
\cdots \wedge \widehat{m_j} \wedge \cdots \wedge \widehat{m_k} \wedge
\cdots \wedge m_d,\quad 0 \le j < k \le d
\end{array}\]
We will also set $\omega_{jk} = -\omega_{kj}$ when $j > k$.  The signs
in the definitions of $\omega_j$ and $\omega_{jk}$ were chosen to make
interior products easy to calculate.  In fact, one can check that if
$m_1^*,\dots,m_d^*$ is the dual basis to $m_1,\dots,m_d$, then
\[ \langle m_k^*,\omega_j\rangle = \cases{0 & $k = j$ \cr \omega_{kj}
& $k \ne j$.\cr}\]
Since $e_i = \sum_{k=1}^d \langle m_k,e_i\rangle m_k^*$, it follows
that
\begin{equation}
 \langle e_i,\omega_j\rangle = \sum_{k \ne j} \langle m_k, e_i\rangle
\omega_{kj}.
\label{eq.interior}
\end{equation}
For later purposes, we also note that as in the proof of Proposition
\ref{prop.dmodule}, the form $\omega_j$ can be written inside
$S\otimes\Lambda^{d-1}{\bf Z}^n$ as
\begin{equation}
\sum_{|J| = d-1}
(-1)^{j+1} \det(e^j_J)(\textstyle{\prod_{l \in J}} z_l^{-1}) dz_J
\label{eq.omegaj}
\end{equation}
where $\det(e_J^j)$ is the determinant of the $(d-1)\times(d-1)$
submatrix of $(\langle m_k,e_l\rangle_{\scriptscriptstyle{1 \le k \le
d, l \in J}})$ obtained by deleting the $j^{\rm th}$ row.

	The next observation is that $S$ has a ${\bf Z}^n$ grading
coming from the action of $({\bf C}^*)^n$ on ${\bf A}^n$.  The graded
pieces are 1-dimensional, each spanned by a single monomial.  This
induces a grading on $S\otimes\Lambda^{d-1}M$ and $S/z_i
S\otimes\Lambda^{d-2}M$, and the map $\gamma$ of (\ref{seq.dminus1})
is a graded homomorphism.  Thus the entire exact sequence has a ${\bf
Z}^n$ grading.  It following that in studying
$\widehat{\Omega}^{d-1}_S \subset S\otimes\Lambda^{d-1}M$, we can
restrict our attention to ``homogeneous'' elements of
$S\otimes\Lambda^{d-1}M$, i.e., elements of the form $\sum_{j=1}^d c_j
z^D\otimes \omega_j$, where $\omega_j$ is as above, $c_j \in {\bf C}$,
and $z^D$ is a monomial.

	Take such an element $\sum_{j=1}^d c_jz^D\otimes\omega_j \in
S\otimes\Lambda^{d-1} M$.  Then, in the exact sequence
(\ref{seq.dminus1}), equation (\ref{eq.interior}) shows that the
$i^{\rm th}$ component of $\gamma(\sum_{j=1}^d c_jz^D\otimes\omega_j)$
is given by
\[ \sum_{j=1}^d\sum_{k\ne j} (c_jz^D \langle m_k,e_i\rangle \bmod
z_i)\otimes \omega_{kj} = \sum_{j<k} (\langle c_jm_k-c_km_j,e_i\rangle
z^D)\bmod z_i)\otimes \omega_{kj}.\]
Thus it follows that $\sum_{j=1}^d c_j z^D\otimes\omega_j$ lies in
$\widehat{\Omega}^{d-1}_S$ if and only if
\begin{equation}
\hbox{$z_i$ divides $\langle c_j m_k-c_km_j,e_i\rangle z^D$ for all
$i$, $j$ and $k$.}
\label{criterion}
\end{equation}

	For $1 \le i \le n$, the above criterion shows that
$\sum_{j=1}^d \langle m_j,e_i\rangle \widehat{z}_i\, \omega_j$ lies in
$\widehat{\Omega}^{d-1}_S$, where $\widehat{z}_i = \prod_{l \ne
i}z_l$.  Furthermore, using equation (\ref{eq.omegaj}), we see that in
$S\otimes\Lambda^{d-1}{\bf Z}^n$, this form equals
\[ \sum_{j=1}^d \langle m_j,e_i\rangle \widehat{z}_i
\Bigl(\textstyle{\sum_{|J| = d-1}} (-1)^{j+1}
\det(e^j_J)(\textstyle{\prod_{l \in J}} z_l^{-1})\Bigr) dz_J, \]
which can be written as
\[ \sum_{|J| = d-1} \Bigl(\textstyle{\sum_{j=1}^d} (-1)^{j+1}\langle
m_j,e_i\rangle \det(e_J^j)\Bigr) \widehat{z}_{J\cup\{i\}} dz_J =
\sum_{|J| = d-1} \det(e_{J\cup\{i\}}) \widehat{z}_{J\cup\{i\}} dz_J =
\Omega_i \]
since $\langle m_j,e_i\rangle$ for $1 \le j \le d$ gives the first
column of the matrix in $\det(e_{J\cup\{i\}})$.  Similarly, if $(i,j)
\in {\cal I}_1$, one can show that $\sum_{j=1}^d \langle
m_j,e_i\rangle \widehat{z}_{ij}\, \omega_j$ lies in
$\widehat{\Omega}^{d-1}_S$ and equals $\Omega_{ij}$.

	Thus $\Omega_i$, $1 \le i \le n$, and $\Omega_{ij}$, $(i,j)
\in {\cal I}_1$ lie in $\widehat{\Omega}^{d-1}_S$.    To prove that
these forms generate $\widehat{\Omega}^{d-1}_S$, suppose $\omega =
\sum_{j=1}^d c_j z^D\otimes \omega_j$ satisfies equation
(\ref{criterion}).  We can assume that at least one $c_j \ne 0$.
There are three cases to consider:

{\sc Case I.} Suppose that $z_i$ divides $z^D$ for all $i$.  It
suffices to consider $\omega = z_1\dots z_n\, \omega_j$, and we can
assume that $e_1,\dots, e_d$ are linearly independent.  Then for $1
\le i \le d$, the $z_i\Omega_i = \sum_{j=1}^n \langle m_j,e_i\rangle
z_1\cdots z_n\, \omega_j$ have the same span as the $z_1\dots z_n\,
\omega_j$ since $(\langle m_j,e_i\rangle_{\scriptscriptstyle{1 \le k,i
\le d}})$ is invertible.  Thus $z_1\cdots z_n\, \omega_j$ lies in the
submodule generated by the $\Omega_i$.

{\sc Case II.} Suppose that $z_i$ doesn't divide $z^D$ for some $i \in
{\cal I}_0$ .  Then (\ref{criterion}) implies that $\langle c_jm_k -
c_km_j, e_i\rangle = 0$ for all $j,k$.  Thus $e_i$ is orthogonal to
the codimension 1 sublattice spanned by $c_jm_k-c_km_j$.  Since $e_i$
is primitive, we see that $e_i$ is determined uniquely up to $\pm$.
In particular, if there were some $j \ne i$ where $z_j$ also didn't
divide $z^D$, then we would have $e_j = -e_i$, which is impossible
since $i \in {\cal I}_0$.  This shows that $\widehat{z}_i$ divides
$z^D$.  Furthermore, $0 = \langle c_jm_k-c_km_j,e_i\rangle =
c_j\langle m_k,e_i\rangle - c_k\langle m_j,e_i\rangle$ shows that for
some constant $\lambda$, we have $c_j = \lambda\langle m_j,e_i\rangle$
for all $j$.  Thus $\omega$ is a multiple of $\Omega_i$.

{\sc Case III.} Suppose that $z_i$ doesn't divide $z^D$ for some $i
\notin {\cal I}_0$.  We can assume $(i,j) \in {\cal I}_1$, and the
argument of Case II shows that $\widehat{z}_{ij}$ divides $z^D$.
Then, as in Case II, we see that $\omega$ is a multiple of
$\Omega_{ij}$.

	Since $\Omega_i = z_j\Omega_{ij}$ and $\Omega_j =
-z_i\Omega_{ij}$ for $(i,j) \in {\cal I}_1$, we only need the the
$\Omega_i$ for $i \in {\cal I}_0$ along with the $\Omega_{ij}$ to
generate $\widehat{\Omega}^{d-1}_S$.  This proves the proposition.
\hfill$\Box$
\medskip

	We get the following description of $(d-1)$-forms on ${\bf P}$
with poles on $X$.

\begin{theo}
\label{theo.dminus1}
Let $X \subset {\bf P}$ be a Cartier divisor defined by
$f = 0$, where $f \in S_\beta$.  Then
\[ H^0({\bf P},\Omega^{d-1}_{\bf P}(X)) =
\Bigl\{ {\sum_{i \in {\cal I}_0} A_i\Omega_i + \sum_{(i,j) \in {\cal
I}_1} A_{ij} \Omega_{ij}\over f} : A_i \in S_{\beta-\beta_0+\beta_i},\
A_{ij} \in S_{\beta-\beta_0+\beta_i + \beta_j} \Bigr\},\]
where $\beta_0 = \sum_{i=1}^n \beta_i$ and $z_i \in S_{\beta_i}$.
However, if $S_\beta \subset B(\Sigma) = \langle
\widehat{z}_\sigma : \sigma \in \Sigma\rangle$ (see Definition
\ref{def.B}), then
\[ H^0({\bf P},\Omega^{d-1}_{\bf P}(X)) =
\Bigl\{ {\sum_{i=1}^n A_i\Omega_i \over f} : A_i \in
S_{\beta-\beta_0+\beta_i} \Bigr\}.\]
\end{theo}

\proof Since $\Omega_i$ has degree $\beta_0 - \beta_i$ and
$\Omega_{ij}$ has degree $\beta_0 - \beta_i - \beta_j$, the first part
of the theorem follows immediately from Propositions
\ref{prop.pform} and \ref{prop.dminus1}.  For the second part,
suppose that $S_\beta \subset B(\Sigma)$, and consider a form $A\,
\Omega_{ij}$ where $A \in S_{\beta -\beta_0 + \beta_i + \beta_j}$.  If
$z^D$ is a monomial that appears in $A$, then $z^D\,\widehat{z}_{ij}
\in S_\beta \subset B(\Sigma)$.  This implies that
$z^D\,\widehat{z}_{ij}$ is divisible by $\widehat{z}_\sigma$ for some
$\sigma \in \Sigma$.  But $\sigma$ can't contain both $e_i$ and $e_j$
since $e_i = -e_j$.  Thus $z_i$ or $z_j$ divides $\widehat{z}_\sigma$,
so that $z_i$ or $z_j$ divides $z^D\,\widehat{z}_{ij}$.  It follows
that $z^D$ is divisible by $z_i$ or $z_j$, and thus $z^D\,\Omega_{ij}$
is a multiple of either $z_i\,\Omega_{ij} = -\Omega_j$ or
$z_j\,\Omega_{ij} = \Omega_j$.  Hence $A\,\Omega_{ij}$ is in the
submodule generated by the $\Omega_i$, and the theorem is proved.
\hfill $\Box$

\begin{rem} {\rm The reader can check that when ${\bf P} = {\bf P}^d$,
the description of $H^0({\bf P},\Omega^{d-1}_{\bf P}(X))$ given above
generalizes equation (4.4) from \cite{griff1}.}
\end{rem}

	To exploit the second part of Theorem \ref{theo.dminus1}, we
need to know when $S_\beta \subset B(\Sigma)$.

\begin{lem}
\label{lem.finB}
 If $\beta \in Cl(\Sigma)$ is the class of an ample divisor
on ${\bf P}$, then $S_\beta \subset B(\Sigma)$.
\end{lem}

\proof Given $z^D = \prod_{i=1}^n z_i^{a_i} \in S_\beta$, our
hypothesis implies that $D = \sum_{i=1}^n a_i D_i$ is an ample divisor on
${\bf P}$.  Thus ${\cal L} = {\cal O}_{\bf P}(D)$ is an ample ${\bf
T}$-linearized invertible sheaf.  Let $\Delta = \Delta({\cal L})
\subset M_{\bf R}$ be its associated convex polytope.  Since $\Delta$
is defined by the inequalities $\langle m,e_i\rangle \ge -a_i$, we
have $0 \in \Delta$ since $a_i \ge 0$.  Then let $\Delta_0$ be the
minimal face of $\Delta$ containing $0$.

	We know that for some $\sigma \in \Sigma$, the face
$\Delta_\sigma$ of $\Delta$ corresponding to $\sigma$ is $\Delta_0$.
We claim that $\widehat{z}_\sigma$ divides our monomial $z^D$.  To see
why this is true, suppose we have some $z_i$ which doesn't divide
$z^D$.  This means $a_i = 0$.  If $\rho_i$ is the 1-dimensional cone
of $\Sigma$ generated by $e_i$, then the corresponding facet of
$\Delta$ is $\Delta_{\rho_i} = \Delta\cap\{m : \langle m,e_i\rangle
\ge 0\}$ since $a_i = 0$.  Thus $0 \in \Delta_{\rho_i}$, which implies
$\Delta_\sigma \subset \Delta_{\rho_i}$ by the minimality of
$\Delta_\sigma$.  It follows that $\rho_i \subset \sigma$.  We have
thus proved that $a_i = 0$ implies $e_i \in \sigma$, and it follows
immediately that $\widehat{z}_\sigma = \prod_{e_i \notin \sigma} z_i$
divides $z^D$.  Thus $z^D \in B(\Sigma)$, and the lemma is proved.
\hfill $\Box$
\medskip

	If we combine this lemma with Theorem \ref{theo.dminus1}, we
get the following useful corollary.

\begin{coro}
\label{coro.dminus1}
Let $X \subset {\bf P}$ be an ample Cartier divisor defined by
$f = 0$, where $f \in S_\beta$.  Then
\[ H^0({\bf P},\Omega^{d-1}_{\bf P}(X)) =
\Bigl\{ {\sum_{i=1}^n A_i\Omega_i \over f} : A_i \in
S_{\beta-\beta_0+\beta_i} \Bigr\}.\]
where $\beta_0 = \sum_{i=1}^n \beta_i$ and $z_i \in S_{\beta_i}$.
\end{coro}

\section{Cohomology of the complement of an ample $\;\;\;\;\;$ divisor}

In this section, ${\bf P}$ will be a $d$-dimensional complete
simplicial toric variety and $X \subset {\bf P}$ will be the zero
locus of a global section of a ${\bf T}$-linearized ample invertible
sheaf ${\cal L}$.  If $\beta \in Cl(\Sigma)$ is the class of ${\cal
L}$, then Lemma \ref{lem.iso} shows that $X$ is defined by an equation
$f = 0$ for some $f \in S_\beta$.  We will also assume that $X$ is
quasi-smooth (by Proposition \ref{prop.generic}, this is true for
generic $f \in S_\beta$).  Our goal is to compute the cohomology of
${\bf P} \backslash X$ in terms of $f \in S$.  We will also study the
cohomology of $X$.  Our results will generalize classical results of
Griffiths, Dolgachev and Steenbrink (see \cite{griff1,dolgach,steen}).

	Since $X$ is a $V$-submanifold of the $V$-manifold ${\bf P}$,
we can compute $H^\cdot({\bf P} \backslash X)$ using the complex
$\Omega_{\bf P}^\cdot({\log X})$ (we always use cohomology with
coefficients in ${\bf C}$).  Furthermore, the Hodge filtration
$F^\cdot$ on $H^{p+q}({\bf P}\backslash X)$ comes from the spectral
sequence $H^q({\bf P},\Omega_{\bf P}^p({\log X})) \Rightarrow
H^{p+q}({\bf P}\backslash X)$ which degenerates at $E_1$ (see \S15 of
\cite{dan1}).  Thus we obtain isomorphisms
\[ Gr_F^p\, H^d({\bf P}\backslash X) \cong H^{d-p}({\bf P},\Omega_{\bf
P}^p(\log X))\]
for $p = 0,\dots,d$.  Morever, $\Omega_{\bf P}^p(\log X)$ has
the following resolution.

\begin{prop}
If $X$ is a quasi-smooth hypersurface of a complete simplicial toric
variety ${\bf P}$, then there is a canonical exact sequence
\[ 0 \rightarrow \Omega^p_{\bf P} ({\rm log}\; X) \rightarrow
\Omega^p_{\bf P} (X) \stackrel{d}{\rightarrow} \Omega^{p+1}_{\bf P} (2X)/
\Omega^{p+1}_{\bf P} (X) \stackrel{d}{\rightarrow} \ldots \]
\[ \ldots \stackrel{d}{\rightarrow} \Omega^{d-1}_{\bf P} ((d -p)X)/
\Omega^{d-1}_{\bf P} ((d-p-1)X)
\stackrel{d}{\rightarrow} \Omega^{d}_{\bf P} ((d -p +1)X)/
\Omega^d_{\bf P} ((d-p)X)
 \rightarrow 0. \]
\end{prop}

\proof The proof of this is similar to the proof of Theorem 6.2 in
\cite{bat.var}. \hfill$\Box$
\medskip

	If we combine the above exact sequence with the
Bott-Steenbrink-Danilov vanishing theorem (see Theorem
\ref{theo.bott}), then we obtain the following corollary.

\begin{coro}
\label{coro.log}
There are natural isomorphisms
\[ Gr_F^p\, H^d({\bf P}\backslash X) \cong H^{d-p}({\bf P},
\Omega_{\bf P}^p ({\rm log}\, X)) \cong
\frac{H^0({\bf P}, \Omega^{d}_{\bf P} ((d -p +1)X)}
{H^0({\bf P}, \Omega^d_{\bf P} ((d -p)X) +
dH^0({\bf P}, \Omega^{d-1}_{\bf P} ((d -p)X)}. \]
\end{coro}

	Before we can state our next result, we need some definitions.

\begin{opr}
{\rm Let $f \in S_\beta$ be a nonzero polynomial.  Then the {\em
Jacobian ideal $J(f) \subset S$} is the ideal of $S = {\bf
C}[z_1,\dots,z_n]$ generated by the partial derivatives $\partial
f/\partial z_1,\dots,\partial f/\partial z_n$.  Also, the {\em
Jacobian ring $R(f)$} is the quotient ring $S/J(f)$.}
\end{opr}

\begin{rem}
{\rm Since $f \in S_\beta$, we have $\partial f/\partial z_i \in
S_{\beta-\beta_i}$, where $z_i \in S_{\beta_i}$.  Thus $J(f)$ is a
graded ideal of $S$, so that $R(f)$ has a natural grading by the class
group $Cl(\Sigma)$.}
\end{rem}

\begin{lem}
\label{lem.fJ}
If $f \in S_\beta$, where $\beta \ne 0$, then $f \in J(f)$.
\end{lem}

\proof  We will prove this using the Euler formulas from Definition
\ref{euler.def}.  Let $\prod_{i=1}^n z_i^{b_i}$ be a monomial
appearing in $f$.  Then, in $Cl(\Sigma)$, we have $\beta =
[\sum_{i=1}^n b_i D_i]$, where $D_i$ is the divisor in ${\bf P}$
corresponding to $e_i$.  If we can find $\phi_1,\dots,\phi_n \in {\bf
C}$ such that $\sum_{i=1}^n \phi_i e_i = 0$, then Lemma
\ref{euler.lem} tells us that
\[ ({\textstyle{\sum_{i=1}^n \phi_ib_i}})\,f = \sum_{i=1}^n \phi_i z_i
{\partial f \over \partial z_i}.\]
Thus we need to find a relation $\sum_{i=1}^n \phi_i e_i = 0$ such
that $\sum_{i=1}^n \phi_i b_i \ne 0$.  To do this, pick $j$ such that
$b_j > 0$.  By completeness, $-e_j$ must lie in some cone $\sigma \in
\Sigma$.  Since $e_j$ can't be a generator of $\sigma$, we get $-e_j =
\sum_{i \ne j} \phi_i e_i$, where $\phi_i \ge 0$, so that setting
$\phi_j = 1$ gives $\sum_{i=1}^n \phi_i e_i = 0$.  Since $b_j > 0$ and
$b_i \ge 0$ for $i \ne j$, we have $\sum_{i=1}^n \phi_i b_i > 0$.
\hfill $\Box$
\medskip

	We can now state the first main result of this section.

\begin{theo}
\label{theo.main1}
Let ${\bf P}$ be a $d$-dimensional complete simplicial toric variety,
and let $X \subset {\bf P}$ be a quasi-smooth ample hypersurface
defined by $f \in S_\beta$.  If $R(f)$ is the Jacobian ring of $f$,
then there is a canonical isomorphism
\[ Gr_F^p H^d({\bf P}\backslash X) \cong R(f)_{(d-p+1)\beta-\beta_0},\]
where $z_i \in S_{\beta_i}$ and $\beta_0 = \sum_{i=1}^n \beta_i$.
\end{theo}

\proof The arguments are similar to those used in the classical case
(see, for example, \cite{pet.steen}).  By Theorem \ref{theo.dform}, we
have
\[ H^0({\bf P},\Omega^d_{\bf P}((d-p+1)X)) =
\Bigl\{ {A\Omega_0 \over f^{d-p+1}} : A \in
S_{(d-p+1)\beta-\beta_0}\Bigr\},\]
so that the map $\phi(A\Omega_0/f^{d-p+1}) = A$ defines a bijection
\[\phi: H^0({\bf P},\Omega^d_{\bf P}((d-p+1)X)) \cong S_{(d-p+1)\beta
-\beta_0}.\]
By Corollary \ref{coro.log}, it suffices to show that the subspace
\[ H^0({\bf P}, \Omega^d_{\bf P} ((d -p)X) + dH^0({\bf P},
\Omega^{d-1}_{\bf P} ((d -p)X) \subset H^0({\bf P},\Omega^d_{\bf
P}((d-p+1)X)) \]
maps via $\phi$ to $J(f)_{(d-p+1)\beta-\beta_0} \subset
S_{(d-p+1)\beta-\beta_0}$.

	If $p = d$, then the desired result follows immediately since
$H^0({\bf P},\Omega_{\bf P}^d)$, $H^0({\bf P},\Omega_{\bf P}^{d-1})$
and $J(f)_{\beta-\beta_0}$ all vanish (the last because $\partial
f/\partial z_i \in S_{\beta-\beta_i}$).  Now assume $p < d$.  Since
$A\Omega_0/f^{d-p} = fA\Omega_0/f^{d-p+1}$, we see that $\phi(H^0({\bf
P},\Omega_{\bf P}^d((d-p)X))) = fS_{(d-p)\beta-\beta_0}$.  It remains
to see what happens to $dH^0({\bf P},\Omega^{d-1}_{\bf P}((d-p)X))$.
By Corollary \ref{coro.dminus1}, we know that
\[ H^0({\bf P},\Omega^{d-1}_{\bf P}((d-p)X)) =
\Bigl\{ {\sum_{i=1}^n A_i\Omega_i \over f^{d-p}} : A_i \in
S_{(d-p)\beta-\beta_0+\beta_i} \Bigr\}.\]
Since $d-p \ne 0$ and $f \in J(f)$ by Lemma \ref{lem.fJ}, the theorem
now follows easily from the following lemma.

\begin{lem}
If $A \in S_{(d-p)\beta-\beta_0+\beta_i}$, then
\[ d\Bigl({A \Omega_i \over f^{d-p}}\Bigr) = {(f \partial A/\partial
z_i - (d-p) A \partial f/\partial z_i)\Omega_0 \over f^{d-p+1}}.\]
\end{lem}

\proof First note that
\begin{equation}
\label{eq.mess}
d\Bigl({A\Omega_i \over f^{d-p}}\Bigr) = {1\over f^{d-p+1}}(f dA
\wedge \Omega_i + f A d\Omega_i - (d-p)df\wedge \Omega_i).
\end{equation}
This equals $B\Omega_0/f^{d-p+1}$ for some $B \in
S_{(d-p+1)\beta-\beta_0}$, where
\[\Omega_0 = \sum_{|I| = d} \det(e_I) \widehat{z}_I dz_I\]
(see Definition \ref{def.omegao}).  Pick $I_0 \subset \{1,\dots,n\}$
such that $|I_0| = d$, $i \in I_0$ and $\det(e_{I_0}) \ne 0$.  To find
$B$, it suffices to determine the coefficient of $dz_{I_0}$ in
(\ref{eq.mess}).

	From Definition \ref{def.omegai}, we have
\[ \Omega_i = \sum_{\scriptstyle{ |J| = d-1 \atop i \notin J}}
\det(e_{J\cup\{i\}}) \widehat{z}_{J\cup\{i\}} dz_J. \]
Since neither $z_i$ nor $dz_i$ appear in $\Omega_i$, $dz_{I_0}$
doesn't appear in $fAd\Omega_i$.  Furthermore, if we set $J_0 = I_0
\backslash \{i\}$, then the coefficient of $dz_{I_0} = dz_i\wedge
dz_{J_0}$ in $fdA\wedge \Omega_i - (d-p)df \wedge \Omega_i$ is
\[f \textstyle{\partial A \over\partial z_i} dz_i \wedge
\det(e_{J_0\cup\{i\}}) \widehat{z}_{J_0\cup\{i\}} dz_{J_0} -
(d-p)A\textstyle{\partial f \over \partial z_i}dz_i
\wedge \det(e_{J_0\cup\{i\}}) \widehat{z}_{J_0\cup\{i\}} dz_{J_0}\]
\[ =  (f \partial A/\partial z_i - (d-p)A\partial f/\partial
z_i)\det(e_{I_0}) \widehat{z}_{I_0} dz_{I_0}.\]
This shows that $B = f \partial A/\partial z_i - (d-p)A\partial
f/\partial z_i$ and completes the proof of the lemma. \hfill$\Box$
\medskip

	We next study the cohomology of the hypersurface $X$.  Our
first result is a Lefschetz theorem.

\begin{prop}
\label{prop.lef}
Let $X$ be a quasi-smooth hypersurface of a $d$-dimensional complete
simplicial toric variety ${\bf P}$, and suppose that $X$ is defined by
$f \in S_\beta$.  If $f \in B(\Sigma)$ (see Definition \ref{def.B}),
then the natural map $i^* : H^i({\bf P}) \to H^i(X)$ is an isomorphism
for $i < d-1$ and an injection for $i = d-1$.  In particular, this
holds if $X$ is an ample hypersurface.
\end{prop}

\proof In the affine space ${\bf A}^n$, $f \in B(\Sigma)$ implies that
$Z(\Sigma) = {\bf V}(B(\Sigma)) \subset {\bf V}(f)$.  Thus ${\bf A}^n
\backslash {\bf V}(f) \subset {\bf A}^n - Z(\Sigma) = U(\Sigma)$ is
affine.  Since ${\bf P} = U(\Sigma)/{\bf D}(\Sigma)$, it follows that
${\bf P} \backslash X = ({\bf A}^n \backslash {\bf V}(f))/{\bf
D}(\Sigma)$ is affine.  This implies $H^i({\bf P}\backslash X) = 0$
for $i > d$ by Corollary 13.6 of
\cite{dan1}.  Now consider the Gysin sequence
\[ \dots \to H^{i-2}(X) \stackrel{i_!}\to H^i({\bf P}) \to H^i({\bf
P}\backslash X) \to H^{i-1}(X) \stackrel{i_!}\to H^{i+1}({\bf P}) \to
\dots\]
(see the proof of Theorem 3.7 of \cite{dan.hov}).  Since the Gysin map
$i_!$ is dual to $i^*$ under Poincar\'e duality, we see that $i^*$ has
the desired property.  Finally, if $X$ is ample, then Lemma
\ref{lem.finB} implies that $f \in B(\Sigma)$. \hfill$\Box$
\medskip

	This result shows that the ``interesting'' part of the
cohomology of $X$ occurs in dimension $d-1$ and consists of those
classes which don't come from ${\bf P}$.  Hence we get the following
definition.

\begin{opr}
{\rm The {\em primitive cohomology group} $PH^{d-1}(X)$ is defined by
the exact sequence
\[ 0 \to H^{d-1}({\bf P}) \to H^{d-1}(X) \to PH^{d-1}(X) \to 0.\]}
\end{opr}

\begin{rem}
{\rm Since $H^{d-1}({\bf P})$ and $H^{d-1}(X)$ have pure Hodge
structures, $PH^{d-1}(X)$ is also pure.  Its Hodge components are
denoted $PH^{p,d-1-p}(X)$.}
\end{rem}

\begin{prop}
\label{prop.prim}
When $X \subset {\bf P}$ is ample, there is an exact sequence
\[ 0 \to H^{d-2}({\bf P}) \stackrel{\cup[X]}\to H^d({\bf P}) \to
H^d({\bf P}\backslash X) \to PH^{d-1}(X) \to 0\]
where $[X] \in H^2({\bf P})$ is the cohomology class of $X$.
\end{prop}

\proof By (3.3) of \cite{pet.steen}, we have a commutative diagram
\begin{equation}
\label{eq.gysin}
\begin{array}{ccc}
H^i({\bf P}) && \\
\ \downarrow\,\scriptstyle{i^*} & \searrow\!\!{}^{\cup [X]} & \\
H^i(X) & \stackrel{i_!}\longrightarrow\hfill & \!\!H^{i+2}({\bf P})
\end{array}
\end{equation}
When $i = d-1$, $\cup [X]$ is an isomorphism by Hard Lefschetz, and it
follows that $PH^{d-1}$(X) is isomorphic to the kernel of $i_! :
H^{d-1}(X) \to H^{d-1}({\bf P})$.  Then the Gysin sequence gives us an
exact sequence
\[ H^{d-2}(X) \stackrel{i_!}\to H^d({\bf P}) \to H^d({\bf P}\backslash
X) \to PH^{d-1}(X) \to 0.\]
Since $\cup [X] : H^{d-2}({\bf P}) \to H^d({\bf P})$ is injective by
Hard Lefschetz and $i^* : H^{d-2}({\bf P}) \to H^{d-2}(X)$ is an
isomorphism by Proposition \ref{prop.lef}, the desired exact sequence
now follows easily using (\ref{eq.gysin}) with $i = d-2$.\hfill$\Box$
\medskip

	The maps in Proposition \ref{prop.prim} are all morphisms of
mixed Hodge structures (with appropriate shifts).  Since $H^i({\bf
P})$ vanishes for $i$ odd and has only $(p,p)$ classes for $i$ even
(see Theorem 3.11 of \cite{oda}), we get the following corollary of
Proposition \ref{prop.prim}.

\begin{coro}
When $p \ne d/2$, there is a natural isomorphism
\[Gr_F^p H^d({\bf P}\backslash X) \cong PH^{p-1,d-p}(X),\]
and when $p = d/2$, there is an exact sequence
\[ 0 \to H^{d-2}({\bf P}) \stackrel{\cup[X]}\to H^d({\bf P}) \to
Gr_F^{d/2}H^d({\bf P}\backslash X) \to PH^{d/2-1,d/2}(X) \to 0.\]
\end{coro}

	If we combine this with Theorem \ref{theo.main1} and replace
$p$ with $p+1$, then we get the following theorem.

\begin{theo}
\label{theo.main}
Let ${\bf P}$ be a $d$-dimensional complete simplicial toric variety,
and let $X \subset {\bf P}$ be a quasi-smooth ample hypersurface
defined by $f \in S_\beta$.  If $R(f)$ is the Jacobian ring of $f$,
then for $p \ne d/2-1$, there is a canonical isomorphism
\[R(f)_{(d-p)\beta-\beta_0} \cong PH^{p,d-1-p}(X)\]
where $z_i \in S_{\beta_i}$ and $\beta_0 = \sum_{i=1}^n \beta_i$.
For $p = d/2-1$, we have an exact sequence
\[ 0 \to H^{d-2}({\bf P}) \stackrel{\cup[X]}\to H^d({\bf P}) \to
R(f)_{(d/2+1)\beta-\beta_0} \to PH^{d/2-1,d/2}(X) \to 0.\]
\end{theo}

\begin{rem}
\label{rem.main}
{\rm Notice that when $d$ is odd, we always have
$R(f)_{(d-p)\beta-\beta_0} \cong PH^{p,d-1-p}(X)$.  The same
conclusion holds whenever ${\bf P}$ is a toric variety with the
property that $\cup [X] : H^{d-2}({\bf P}) \to H^d({\bf P})$ is an
isomorphism.  The latter holds when ${\bf P}$ is a weighted projective
space, which explains why the classical case is so nice.  Notice also
that in all cases, we always have a surjection
\[R(f)_{(d-p)\beta-\beta_0} \to PH^{p,d-1-p}(X) \to 0.\]}
\end{rem}

\begin{rem}
{\rm One consequence of our results is that there is a natural map
\[ H^d({\bf P}) \to R(f)_{(d/2+1)\beta-\beta_0}.\]
It would be interesting to have an explicit description of this map.}
\end{rem}

\section{Cohomology of affine hypersurfaces}

        In a recent paper \cite{bat.var}, the first author obtained
some results on the cohomology of the affine hypersurface $Y =
X\cap{\bf T} \subset {\bf T}$, where ${\bf T}$ is the torus contained
in a projective  toric variety ${\bf P}$ and $X \subset {\bf
P}$ is an ample hypersurface.  We will show how the results of
\cite{bat.var} can be expressed in terms of various graded ideals of
the ring $S = {\bf C}[z_1,\dots,z_n]$.  We begin with a definition.

\begin{opr} {\rm If $X \subset {\bf P}$ is a hypersurface, then we let
$Y = X\cap {\bf T} \subset {\bf T}$ be its intersection with the torus
${\bf T} \subset {\bf P}$.  The primitive cohomology group
$PH^{d-1}(Y)$ is then defined to be the cokernel of the map
$H^{d-1}({\bf T}) \to H^{d-1}(Y)$, where as usual we use cohomology
with coefficients in ${\bf C}$.}
\end{opr}

\begin{rem} {\rm In \cite{dan.hov}, it is shown that $H^{d-1}({\bf T})
\to H^{d-1}(Y)$ is injective.  Note also that $PH^{d-1}(Y)$ has a
natural mixed Hodge structure.  The Hodge filtration on $PH^{d-1}(Y)$
will be denoted $F^\cdot$.}
\end{rem}

        One of the main results of \cite{bat.var} is a description of
the Hodge filtration on $PH^{d-1}(Y)$.  To state this in terms of the
ring $S$, we need the following ideal of $S$.

\begin{opr} {\rm Given $f \in S_\beta$, let $J_0(f) \subset S$ denote
the ideal generated by $z_i \partial f/\partial z_i$ for $1 \le i \le
n$.  We then let $R_0(f)$ denote the quotient ring $S/J_0(f)$.}
\end{opr}

\begin{rem} {\rm Since $f \in S_\beta$, we see that $J_0(f)$ is a
graded ideal of $S$, and hence $R_0(f)$ has a natural grading by the
class group $Cl(\Sigma)$.}
\end{rem}

\begin{theo}
\label{theo.aff}
{\rm (\cite{bat.var})} If $X \subset {\bf P}$ is a nondegenerate (see
Definition \ref{def.nondeg}) ample divisor defined by $f \in S_\beta$
and $Y = X\cap {\bf T}$ is the corresponding affine hypersurface in
the torus ${\bf T}$, then there is a natural isomorphism
\[ Gr^p_F PH^{d-1}(Y) \cong R_0(f)_{(d-p)\beta}.\]
\end{theo}

\proof We first show how certain constructions in \cite{bat.var} can
be formulated in terms of the ring $S$.  Since $X \subset {\bf P}$ is
ample, we get the $d$-dimensional convex polyhedron $\Delta \subset
M_{\bf R}$.  Recall that $\Delta$ is defined by the inequalities
$\langle m,e_i\rangle \ge -b_i$, where $b_i \ge 0$ and $X$ is linearly
equivalent to $D = \sum_{i=1}^n b_i D_i$.

        As in \cite{bat.var}, we define the ring $S_\Delta$ to be the
subring of ${\bf C}[t_0,t_1^{\pm1},\dots,t_d^{\pm1}]$ spanned over
${\bf C}$ by all Laurent monomials of the form $t_0^k t^m =t_0^k
t_1^{m_1}\cdots t_d^{m_d}$ where $k \ge 0$ and $m \in k\Delta$.  This
ring is graded by setting $\deg(t_0^kt^m) = k$.  We should think of
$t_1,\dots,t_d$ as coordinates on the torus ${\bf T}$ and $t_0$ as an
auxillary variable to keep track of the grading.

        A first observation is that there is a natural isomorphism of
graded rings
\begin{equation}
\rho :  S_\Delta \cong \bigoplus_{k=0}^\infty S_{k\beta} \subset S
\label{eq.translate}
\end{equation}
which is defined by
\[ \rho(t_0^k t^m) = \prod_{i=1}^n z_i^{kb_i + \langle m,e_i\rangle}.\]
Since the integer points of $k\Delta$ naturally give a basis of
$H^0({\bf P},{\cal O}_{\bf P}(kD))$, Proposition 1.1 of \cite{cox} shows
that we get the desired ring isomorphism.

        In particular, $f \in S_\beta$ corresponds to an element
$\sum_{m \in \Delta\cap M} a_m t_0 t^m \in (S_\Delta)_1$.  Setting
$t_0 = 1$, we get the Laurent polynomial $g = \sum_{m \in \Delta\cap
M} a_m t^m$ formed using the lattice points of $\Delta$.  Since
$t_1,\dots,t_d$ are coordinates on the torus ${\bf T}$, one can show
that $Y$ is naturally isomorphic to the subvariety of ${\bf T}$
defined by $g = 0$.

        Following \cite{bat.var}, we use the Laurent polynomial $g$ to
define $F = t_0 g -1$, and then we set $F_i = t_i \partial
F/\partial t_i$ for $i = 0,1,\dots,d$.  Note that $F_i \in
(S_\Delta)_1$ for all $i$.  Finally, we define $J_{g,\Delta}$ to be the
ideal of $S_\Delta$ generated by the $F_i$.  Then Corollary 6.10 of
\cite{bat.var} gives an isomorphism
\begin{equation}
Gr_F^p PH^{d-1}(Y) \cong (S_\Delta/J_{g,\Delta})_{d-p}.
\end{equation}
To prove the theorem, it suffices to show that for all $k \ge 0$, the
isomorphism $\rho$ of (\ref{eq.translate}) maps the graded piece
$(J_{g,\Delta})_k \subset (S_\Delta)_k$ to $J_0(f)_{k\beta} \subset
S_{k\beta}$.

        First observe that $F_0 = t_0 g$ maps to $f$ under the map
$\rho$.  To write this explicitly, we introduce the following
notation: given $m \in \Delta\cap M$, let $z^{D(m)} = \prod_{i=1}^n
z_i^{b_i + \langle m,e_i\rangle}$.  Thus $f = \sum_{m \in \Delta\cap
M} a_m z^{D(m)}$.  Next note that $t_1,\dots,t_d$ are a basis of the
character group of ${\bf T}$ and hence determine a basis of $M$.  If
$h_1,\dots,h_d \in N$ is the dual basis, then one computes that for $i
> 0$, the polynomial $F_i = t_0t_i \partial g/\partial t_i$ maps to
\[ \tilde f_i = \sum_{m\in \Delta\cap M} a_m\langle m,h_i\rangle
z^{D(m)}.\]
In contrast, note that
\[ z_i \partial f/\partial z_i = \sum_{m\in\Delta\cap M} a_m(b_i +
\langle m,e_i\rangle) z^{D(m)} = b_i f + \sum_{m\in\Delta\cap M} a_m
\langle m,e_i\rangle z^{D(m)}.\]
It suffices to show that each $z_i\partial f/\partial z_i$ is a
${\bf C}$-linear combination of $f,\tilde f_1,\dots,\tilde f_d$ and
conversely.

        To prove this, first note that $e_i$ can be expressed in terms
of the $h_1,\dots,h_d$, and hence $z_i\partial f/\partial z_i$ is a
${\bf C}$-linear combination of $f,\tilde f_1,\dots,\tilde f_d$.
Going the other way, we can label $e_1,\dots,e_n$ so that
$e_1,\dots,e_d$ are linearly independent.  Then $h_i$ can be expressed
in terms of $e_1,\dots,e_d$, and it follows easily that $f,\tilde
f_1,\dots,\tilde f_d$ are ${\bf C}$-linear combinations of
$f,z_1\partial f/\partial z_1,\dots,z_n\partial f/\partial z_n$.  To
complete the proof, observe that $f \in J_0(f)$ follows easily from
the proof of Lemma \ref{lem.fJ}.\hfill $\Box$
\medskip

	We next study the weight filtration on $PH^{d-1}(Y)$.  Since
$Y$ is quasi-smooth, we have $W_{d-2}H^{d-1}(Y) = 0$, which implies
that $W_{d-2}PH^{d-1}(Y) = 0$.  It follows that $W_{d-1}PH^{d-1}(Y)$
has a pure Hodge structure.  We can identify this Hodge structure as
follows.

\begin{prop}
\label{prop.weight}
If $X \subset {\bf P}$ is a nondegenerate ample hypersurface, then
there is a natural isomorphism of Hodge structures
\[ PH^{d-1}(X) \cong W_{d-1}PH^{d-1}(Y).\]
\end{prop}

\proof Let $D = {\bf P} \setminus {\bf T}$, and recall from Theorem
\ref{poin} that the complex $\Omega_{\bf P}^\cdot(\log D)$ has a weight
filtration $W_\cdot$ with the property that
\[ Gr^W_k \Omega^\cdot_{\bf P}(\log D) \cong \bigoplus_{\dim \tau = k}
\Omega_{{\bf P}_\tau}^\cdot.\]
The spectral sequence of this filtered complex gives the weight
filtration on $H^\cdot({\bf T})$, and since the spectral sequence
degenerates at $E_2$, we get an exact sequence
\[ \bigoplus_{i=1}^n H^{d-3}(D_i) \to H^{d-1}({\bf P}) \to Gr_{d-1}^W
H^{d-1}({\bf T}) \to 0,\]
where $D = \sum_{i=1}^n D_i$.  Notice also that $Gr_{d-1}^W
H^{d-1}({\bf T}) = W_{d-1}H^{d-1}({\bf T})$ since ${\bf T}$ is smooth.
Since $X \subset {\bf P}$ is nondegenerate, the same argument applies
to $D\cap X = X \setminus Y$, so that we also have an exact sequence
\[ \bigoplus_{i=1}^n H^{d-3}(D_i\cap X) \to H^{d-1}(X) \to Gr_{d-1}^W
H^{d-1}(Y) \to 0,\]
and as noticed earlier, $Gr_{d-1}^W H^{d-1}(Y) = W_{d-1} H^{d-1}(Y)$.

	To see how this applies to primitive cohomology, consider the
commutative diagram:
\[
\begin{array}{ccccccc}
&& 0 && 0 && \\
&& \uparrow && \uparrow && \\
&& PH^{d-1}(X) & \stackrel{\alpha}\to & W_{d-1}PH^{d-1}(Y) && \\
&& \uparrow && \uparrow && \\
\bigoplus_{i=1}^n H^{d-3}(D_i\cap X) & \to & H^{d-1}(X) & \to &
W_{d-1}H^{d-1}(Y) & \to & 0\\
\uparrow && \uparrow && \uparrow && \\
\bigoplus_{i=1}^n H^{d-3}(D_i) & \to & H^{d-1}({\bf P}) & \to &
W_{d-1}H^{d-1}({\bf T}) & \to & 0
\end{array}
\]
The columns are exact by the definition of primitive cohomology and
the strictness of the weight filtration, and we've already seen that
the bottom two rows are exact.  It follows easily that the map $\alpha
: PH^{d-1}(X) \to W_{d-1}PH^{d-1}(Y)$ exists and is surjective.
Notice also that $\alpha$ is a morphism of Hodge structures.

	However, each $D_i$ is a $(d-1)$-dimensional complete
simplicial toric variety, and $D_i\cap X \subset D_i$ is quasi-smooth
since $X$ is nondegenerate.  Thus Proposition \ref{prop.lef} implies
that
\[ \bigoplus_{i=1}^n H^{d-3}(D_i) \to \bigoplus_{i=1}^n
H^{d-3}(D_i\cap X)\]
is an isomorphism.  An easy diagram chase then shows that $\alpha$ is
injective, and the proposition is proved.\hfill$\Box$
\medskip

	In order to interpret this proposition in terms of the
polynomial ring $S$, we will need the following ideal.

\begin{opr} {\rm Given the ideal $J_0(f) = \langle z_1\partial
f\partial z_1,\dots,z_n \partial f/\partial z_n\rangle \subset S$, we
get the ideal quotient $J_1(f) = J_0(f)\colon z_1\cdots z_n$.  We
put $R_1(f) = S/J_1(f)$. }
\end{opr}

\begin{theo} {\rm (\cite{bat.var})} If $X \subset {\bf P}$ is a
nondegenerate ample divisor defined by $f \in S_\beta$, then there is
a natural isomorphism
\[ PH^{p,d-1-p}(X) \cong R_1(f)_{(d-p)\beta-\beta_0},\]
where $z_i \in S_{\beta_i}$ and $\beta_0 = \sum_{i=1}^n \beta_i$.
\label{theo.r1}
\end{theo}

\proof As in Definition 2.8 of \cite{bat.var}, consider the ideal
$I^{(1)}_\Delta \subset S_\Delta$ spanned by all monomials $t_0^k t^m$
such that $m$ is an interior point of $k\Delta$.  Then let $H = \oplus
H_i$ denote the image of the homogeneous ideal $I^{(1)}_\Delta$ in the
graded Artinian ring $S_\Delta/J_{g,\Delta}$, where $J_{g,\Delta}$ is
as in the proof of Theorem \ref{theo.aff}.  Proposition 9.2 of
\cite{bat.var} tells us that under the isomorphism
$(S_\Delta/J_{f,\Delta})_{d-p} \cong Gr_F^p PH^{d-1}(Y)$, the subspace
$H_{d-p}$ maps to $W_{d-1} Gr_F^p PH^{d-1}(Y)$.  If we combine this
with Proposition \ref{prop.weight}, then we get an isomorphism
\[ H_{d-p} \cong PH^{p,d-1-p}(X). \]

        Under the isomorphism $\rho$ of (\ref{eq.translate}), a
monomial $t_0^k t^m$ maps to $z^{D(m)} = \prod_{i=1}^n z_i^{b_i +
\langle m,e_i\rangle}$.  Since $m$ is in the interior of $k\Delta$ if
and only if $\langle m,e_i\rangle > -b_i$ for all $i$, we see that
$t_0^kt^m \in (I^{(1)}_\Delta)_k$ exactly when $z^{D(m)}$ is divisible
by $z_1\cdots z_n$.  Thus $\rho((I^{(1)}_\Delta)_k) = \langle z_1\cdots
z_n\rangle_{k\beta}$.  Since we know by the proof of Theorem
\ref{theo.aff} that $\rho$ maps $(J_{g,\Delta})_k$ to
$J_0(f)_{k\beta}$, it follows easily that $H_k \subset
(S_\Delta/J_{g,\Delta})_k$ is isomorphic to the image of $\langle
z_1\cdots z_n\rangle_{k\beta}$ in $(S/J_0(f))_{k\beta}$.  This last
subspace is isomorphic to $(S/J_1(f))_{k\beta-\beta_0} =
R_1(f)_{k\beta-\beta_0}$, and the theorem is proved.  \hfill$\Box$

\begin{rem} {\rm It is interesting to compare this result to Theorem
\ref{theo.main}, which gives a natural surjection
\begin{equation}
 R(f)_{(d-p)\beta-\beta_0} \to PH^{p,d-1-p}(X) \to 0
\label{eq.main}
\end{equation}
when $X$ is quasi-smooth.  The theorem just proved shows that,
under the stronger hypothesis that $X$ is nondegenerate, there is an
isomorphism
\[ PH^{p,d-1-p}(X) \cong R_1(f)_{(d-p)\beta-\beta_0}.\]
One can show that the composition of these maps in induced by the
obvious inclusion of ideals
\[ J(f) = \langle \partial f/\partial z_i\rangle \subset \langle
z_i \partial f/\partial z_i \rangle\colon z_1\cdots z_n = J_1(f).\]
Since the map of (\ref{eq.main}) is an isomorphism for $p \ne
d/2-1$, it follows that the ideals $J(f)$ and $J_1(f)$ agree in
degrees $(d-p)\beta-\beta_0$ for $p \ne d/2-1$, though for $p =
d/2-1$, we get an exact sequence
\[ 0 \to H^{d-2}({\bf P}) \stackrel{\cup [X]}\to H^d({\bf P}) \to
(J_1(f)/J(f))_{(d/2+1)\beta-\beta_0} \to 0.\]

	When ${\bf P}$ is a weighted projective space, the ideals
$J(f)$ and $J_1(f)$ are equal in all degrees.  For in this case, $f$
being quasi-smooth means that the $\partial f/\partial z_i$ form a
regular sequence, while $f$ being nondegenerate means that the $z_i
\partial f/\partial z_i$ form a regular sequence.  Then standard
results from commutative algebra (see part ($\gamma$) of (1.2) of
\cite{scheja}) imply that $\langle \partial f/\partial z_i\rangle =
\langle z_i\partial f/\partial z_i\rangle\colon z_1\cdots z_n$, which
gives the desired equality.  In the general case, the precise relation
between $J(f)$ and $J_1(f)$ is not well understood.}
\end{rem}

\section{A generalized Euler short exact sequence and $\;\;\;\;\;$
applications}

We begin with a generalization of the classical Euler short exact
sequence.

\begin{theo} Let $D_1, \ldots, D_n$ the irreducible components of $D =
{\bf P} \setminus {\bf T}$, and let $d$ be the dimension of ${\bf P}$.
Then there exists the following short exact sequence
\begin{equation}
 0 \rightarrow \Omega^1_{\bf P} \rightarrow \bigoplus_{i =1}^n
{\cal O}_{\bf P}(-D_i) \rightarrow {\cal O}^{n-d}_{\bf P}
\rightarrow 0.
\label{euler.seq}
\end{equation}
\end{theo}

\begin{rem}
{\rm When ${\bf P}$ is projective space, the above short exact
sequence coincides with the well-known Euler short exact sequence.  So
we call (\ref{euler.seq}) the {\em generalized Euler exact sequence}.}
\end{rem}

\proof There are the following two short exact sequences:
\begin{equation}
 0 \rightarrow \Omega^1_{\bf P} \rightarrow \Omega^1_{\bf P}(\log D)
 \stackrel{r}{\rightarrow} \bigoplus_{i =1}^n
{\cal O}_{D_i} \rightarrow  0
\label{log.seq}
\end{equation}
and
\begin{equation}
 0 \rightarrow \bigoplus_{i =1}^n
{\cal O}_{\bf P}(-D_i) \rightarrow {\cal O}_{\bf P}^{n}
\stackrel{p}{\rightarrow} \bigoplus_{i =1}^n
{\cal O}_{D_i}  \rightarrow 0,
\label{ideal.seq}
\end{equation}
where $r$ is the Poincar\'e residue map.

	The short exact sequence
\[   0 \rightarrow M  \rightarrow {\bf Z}^n \rightarrow Cl(\Sigma)
\rightarrow  0 \]
shows that $Cl(\Sigma)$ has rank $n-d$.  Since $\Omega^1_{\bf P}(\log
D) \cong \Lambda^1 M \otimes_{\bf Z} {\cal O}_{\bf P}$, we can tensor
this sequence with ${\cal O}_{\bf P}$ to obtain
\[  0 \rightarrow \Omega^1_{\bf P}(\log {\bf D})
\stackrel{s}{\rightarrow} {\cal O}_{\bf P}^n  {\rightarrow}
{\cal O}_{\bf P}^{n-d} \rightarrow  0.   \]
By global properties of the Poincar\'e residue map (Section
6), one has $r = p \circ s$. So $s$ induces an injective map $\imath :
\Omega^1_{\bf P} \hookrightarrow \bigoplus_{i =1}^n {\cal O}_{\bf
P}(-D_i)$, and then the short exact sequences (\ref{log.seq}) and
(\ref{ideal.seq}) fit into the following commutative diagram:

\begin{center}

\begin{tabular}{ccccccccc}
  & $   $ & $ 0 $ & $   $ & $ 0 $ & $     $ & $  $ & $ $ &    \\
 & $   $ & $ \downarrow $ & $   $ & $ \downarrow $ & $
 $ & $  $ & $ $ &    \\
$ 0 $ & $ \rightarrow   $ & $ \Omega^1_{\bf P} $ & $
\stackrel{\imath}\rightarrow
$ & $ \bigoplus_{i=1}^n
{\cal O}_{\bf P}(-D_i)  $ & $ \rightarrow
$ & $ {\cal Q} $ & $ \rightarrow $ & $ 0 $ \\
 & $   $ & $ \downarrow $ & $   $ & $ \downarrow $ & $
 $ & $ \downarrow  $ & $ $ &    \\
$ 0 $ & $ \rightarrow   $ & $ \Omega^1_{\bf P}(\log {\bf D}) $ & $
\stackrel{s}{\rightarrow}
$ & $ {\cal O}_{\bf P}^n $ & $ \rightarrow     $ & $ {\cal O}_{\bf P}^{n-d}
$ & $ \rightarrow $ & $ 0$  \\
   & $   $ & $\; \; \downarrow {\scriptstyle r}
   $ & $   $ & $ \;\; \downarrow {\scriptstyle p} $ & $
  $ & $ \downarrow  $ & $ $ &    \\
$ 0 $ & $ \rightarrow   $ & $ \bigoplus_{i=1}^n {\cal O}_{D_i} $ & $
=   $ & $ \bigoplus_{i=1}^n {\cal O}_{D_i}  $ & $ \rightarrow     $ &
$ 0 $ & $  $ &    \\
  & $   $ & $ \downarrow $ & $   $ & $ \downarrow $ & $     $ & $
 $ & $ $ &    \\
   & $   $ & $ 0 $ & $   $ & $ 0 $ & $     $ & $  $ & $ $ &    \\
  \end{tabular}
\end{center}

\noindent Since the first two columns are exact, the snake lemma
implies that the sheaf ${\cal Q}$ is isomorphic to ${\cal O}_{\bf
P}^{n-d}$.
\hfill $\Box$
\medskip

	As a first application of the Euler exact sequence, we will
show how to find generators of $H^0({\bf P},\Omega^{d-1}_{\bf P}(X))$,
where $X \subset {\bf P}$ is an ample hypersurface defined by $f \in
S_\beta$.  If we apply the functor $Hom(*, \Omega^d_{\bf P}(X))$ to
(\ref{euler.seq}), then we obtain the exact sequence
\[  0 \rightarrow (\Omega^d_{\bf P}(X))^{n-d} \rightarrow \bigoplus_{i =1}^n
Hom ({\cal O}_{\bf P}(-D_i),\Omega^d_{\bf P}(X))
\rightarrow \Omega^{d-1}_{\bf P}(X) \rightarrow 0, \]
since $Hom(\Omega^k_{\bf P}, \Omega^d_{\bf P}) = \Omega^{d-k}_{\bf P}$ (see
\cite{dan1}).
Then we get the short exact sequence of global sections
\[  0 \rightarrow H^0({\bf P}, (\Omega^d_{\bf P}(X))^{n-d})
 \rightarrow \bigoplus_{i =1}^n H^0( {\bf P},
 Hom ({\cal O}_{\bf P}(-D_i),\Omega^d_{\bf P}(X)))
\rightarrow H^0({\bf P}, \Omega^{d-1}_{\bf P}(X)) \rightarrow 0 \]
since $H^1({\bf P},\Omega^d_{\bf P}(X)) = 0$ by the
Bott-Steenbrink-Danilov vanishing theorem.

	However, by Remark \ref{rem.omegad}, we know that $\Omega_{\bf
P}^d \cong {\cal O}_{\bf P}(-D)$, and it follows that we have an
isomorphism of ${\bf T}$-linearized sheaves
\[ Hom ({\cal O}(-D_i),\Omega^d_{\bf P}(X)) \cong {\cal O}_{\bf P}(X -
D +D_i).\]
Then Lemma \ref{lem.iso} implies that we have an isomorphism
\[ H^0( {\bf P},  Hom ({\cal O}(-D_i),\Omega^d_{\bf P}(X))) \cong
S_{\beta - \beta_0  + \beta_i}, \]
where as usual $\beta_i$ is the class of $D_i$ and $\beta =
\sum_{i=1}^n \beta_i$ is the class of $D$.  We have thus proved the
following result.

\begin{prop} When $X$ is an ample hypersurface of ${\bf P}$ defined by
$f \in S_\beta$, then there exists a surjective homomorphism
\[ \bigoplus_{i =1}^n S_{\beta - \beta_0  + \beta_i} \rightarrow
H^0({\bf P}, \Omega^{d-1}_{\bf P}(X)). \]
\end{prop}

\begin{rem}
{\rm The reader should compare this proposition with Theorem
\ref{theo.dminus1} and Corollary \ref{coro.dminus1}.}
\end{rem}

We conclude this section with a discussion of the tangent sheaf of a
toric variety.

\begin{opr}
{\rm Let ${\cal T}_{\bf P}$ be the sheaf
$Hom(\Omega^1_{\bf P},{\cal O}_{\bf P})$. We call ${\cal T}_{\bf P}$ the
{\em tangent sheaf} of ${\bf P}$.  }
\end{opr}

\begin{rem}
{\rm If ${\bf P}$ is smooth, then
${\cal T}_{\bf P}$ coincides with the usual tangent sheaf $\Theta_{\bf
P}$ of ${\bf P}$.}
\end{rem}

Applying $Hom(*,{\cal O}_{\bf P})$ to (\ref{euler.seq}), we obtain
the short exact sequence
\[ 0 \rightarrow {\cal O}_{\bf P}^{n-d} \rightarrow
\bigoplus_{i =1}^n {\cal O}(D_i) \rightarrow {\cal T}_{\bf P}
\rightarrow 0. \]
Since $H^1({\bf P}, {\cal O}_{\bf P}) = 0$, we get
\[ {\rm dim}\, H^0({\bf P}, {\cal T}_{\bf P}) =
\sum_{i =1}^n ({\rm dim}\, S_{\beta_i} - 1) + d. \]
It was proved by second author in \cite{cox} that
\[ {\rm dim\,Aut}({\bf P}) =
\sum_{i =1}^n ({\rm dim}\, S_{\beta_i} - 1) + d. \]
Thus the global sections of ${\cal T}_{\bf P}$ can be identified
with the Lie algebra of ${\rm Aut}({\bf P})$.

\section{Moduli of ample hypersurfaces}

	This section will study how the Jacobian ring is related to
the moduli of the hypersurfaces $X \subset {\bf P}$ coming from
sections of an ample invertible sheaf ${\cal L}$ on ${\bf P}$.  As
usual, we assume that ${\bf P}$ is a complete simplicial toric
variety.

	We first study the automorphisms of ${\bf P}$ which preserve a
given divisor class.  Let ${\rm Aut}({\bf P})$ denote the automorphism
group of ${\bf P}$.

\begin{opr}
{\rm Given $\beta \in Cl(\Sigma)$, let ${\rm Aut}_\beta({\bf P})$
denote the subgroup of ${\rm Aut}({\bf P})$ consisting of those
automorphisms which preserve $\beta$.  }
\end{opr}

\begin{rem}
{\rm If ${\rm Aut}{}^0({\bf P})$ is the connected component of the
identity of ${\rm Aut}({\bf P})$, then the results of \S3 of
\cite{cox} imply that ${\rm Aut}{}^0({\bf P})$ is a subgroup of finite
index in ${\rm Aut}_\beta({\bf P})$.}
\end{rem}

	When we describe ${\bf P}$ as the quotient $U(\Sigma)/{\bf
D}(\Sigma)$, note that ${\rm Aut}({\bf P})$ doesn't act on
$U(\Sigma)$.  However, in \cite{cox}, it is shown that there is an
exact sequence
\[ 1 \to {\bf D}(\Sigma) \to \widetilde{\rm Aut}({\bf P}) \to {\rm
Aut}({\bf P}) \to 1\]
where $\widetilde{\rm Aut}({\bf P})$ is the group of automorphisms of
${\bf A}^n$ which preserve $U(\Sigma)$ and normalize ${\bf
D}(\Sigma)$.  An element $\phi \in \widetilde{\rm Aut}({\bf P})$
induces an automorphism $\phi : S \to S$ which for all $\gamma \in
Cl(\Sigma)$ satisfies $\phi(S_\gamma) = S_{\phi(\gamma)}$.

\begin{opr}
{\rm Given $\beta \in Cl(\Sigma)$, let $\widetilde{\rm
Aut}{}_\beta({\bf P})$ denote the subgroup of $\widetilde{\rm
Aut}({\bf P})$ consisting of those automorphisms which preserve
$\beta$.}
\end{opr}

	The group $\widetilde{\rm Aut}_\beta({\bf P})$ has the
following obvious properties.

\begin{lem} There is a canonical exact sequence
\[ 1 \to {\bf D}(\Sigma) \to \widetilde{\rm Aut}_\beta({\bf P}) \to
{\rm Aut}_\beta({\bf P}) \to 1.\]
Furthermore, there is a natural action of $\widetilde{\rm
Aut}_\beta({\bf P})$ on $S_\beta$.
\end{lem}

\begin{rem}
{\rm Let $\widetilde{\rm Aut}{}^0({\bf P})$ be the connected component
of the identity of $\widetilde{\rm Aut}({\bf P})$.  In \cite{cox}, it
is shown that $\widetilde{\rm Aut}{}^0({\bf P})$ is naturally
isomorphic to the group ${\rm Aut}_g(S)$ of $Cl(\Sigma)$-graded
automorphisms of $S$.  Then $\widetilde{\rm Aut}{}^0({\bf P}) \subset
\widetilde{\rm Aut}_\beta({\bf P})$, and the action of $\widetilde{\rm
Aut}_\beta({\bf P})$ on $S_\beta$ is compatible with the action of
${\rm Aut}_g(S)$. }
\end{rem}

	If $\beta \in Cl(\Sigma)$ is an ample class, then we know that
$X = \{p \in {\bf P}: f(p) = 0\} \subset {\bf P}$ is quasi-smooth for
generic $f \in S_\beta$ (see Proposition \ref{prop.generic}).  Then
\[ \{f \in S_\beta : f\ \hbox{is quasi-smooth}\}/\widetilde{\rm
Aut}_\beta({\bf P})\]
should be the coarse moduli space of quasi-smooth hypersurfaces in
${\bf P}$ in the divisor class of $\beta$.  The problem is that
$\widetilde{\rm Aut}_\beta({\bf P})$ need not be a reductive group, so
that the quotient may not exist.  However, it is well-known (see, for
example, \S2 of \cite{cox.tu.don}) that there is a nonempty invariant
open set
\[ U \subset \{f \in S_\beta : f\ \hbox{is quasi-smooth}\}\]
such that the geometric quotient
\[ U/\widetilde{\rm Aut}_\beta({\bf P})\]
exists.

\begin{opr}
{\rm We call the quotient $U/\widetilde{\rm Aut}_\beta({\bf P})$ a
{\em generic coarse moduli space for hypersurfaces of ${\bf P}$ with
divisor class $\beta$}.}
\end{opr}

	We can relate the Jacobian ring $R(f)$ to the generic coarse
moduli space as follows.

\begin{prop}
If $\beta$ is ample and $f \in S_\beta$ is generic, then $R(f)_\beta$
is naturally isomorphic to the tangent space of the generic coarse
moduli space of quasi-smooth hypersurfaces of ${\bf P}$ with divisor
class $\beta$.
\end{prop}

\proof First note that $\widetilde{\rm Aut}{}^0({\bf P}) \subset
\widetilde{\rm Aut}_\beta({\bf P})$ has finite index.  Thus, by
shrinking $U$ if necessary, we may assume that
\[ U/\widetilde{\rm Aut}{}^0({\bf P}) \to U/\widetilde{\rm
Aut}_\beta({\bf P})\]
is \'etale.  Hence it suffices to identify $R(f)_\beta$ with the
tangent space to $U/\widetilde{\rm Aut}{}^0({\bf P})$.  Shrinking $U$
further, we may assume that the map
\[ U \to U/\widetilde{\rm Aut}{}^0({\bf P})\]
is smooth (see \S2 of \cite{cox.tu.don}).  Then the tangent space to
a point of the generic moduli space is naturally isomorphic to the
quotient of $S_\beta$ (= tangent space of $U$) modulo the tangent
space to the orbit of ${\rm Aut}_g(S) = \widetilde{\rm Aut}{}^0({\bf
P})$ acting on $f \in S_\beta$.

	Hence, to prove the proposition, we need tos show that
$J(f)_\beta$ is the tangent space to the orbit of $f$.  But the
tangent space to the orbit is given by the action of the Lie algebra
of ${\rm Aut}_g(S)$ on $f$.  Since the Lie algebra consists of all
derivations of $S$ which preserve the grading, its elements can be
written in the form $\sum_{i=1}^n A_i \partial/\partial z_i$, where
$A_i$ and $z_i$ lie in the same graded piece $S_{\beta_i}$ of $S$ for
all $i$.  Thus the action of the Lie algebra on $f$ gives the subspace
$\{ \sum_{i=1}^n A_i \partial f/\partial z_i : A_i \in S_{\beta_i}\} =
J(f)_\beta$, and the proposition is proved.\hfill$\Box$


\begin{thebibliography}{99}

\bibitem{aud} M. Audin, {\sl The Topology of Torus Actions
on Symplectic Manifolds}, Progress in Math. {\bf 93},
Birkh\"auser, Boston-Basel-Berlin, 1991.

\bibitem{bat.class} V. Batyrev, {\em On the classification of smooth
projective toric varieties}, T\^ohoku Math. J. {\bf 43} (1991),
569--585.

\bibitem{bat.var} V. Batyrev, {\em Variations of the mixed Hodge
structure  of affine hypersurfaces in algebraic tori},
Duke Math. J {\bf 69} (1993), 349--409.

\bibitem{BE} D. Buchsbaum, D. Eisenbud, {\em Algebraic structures for
finite free resolutions, and some structure theorems for ideals of
codimension 3}, Amer. J. Math. {\bf 99} (1977), 447--485.

\bibitem{cox} D. Cox, {\em The homogeneous coordinate ring of a toric
variety}, to appear.

\bibitem{cox.tu.don} D. Cox, R. Donagi, L. Tu, {\em Variational
Torelli implies generic torelli}, Invent. Math. {\bf 88} (1987),
439--446.

\bibitem{dan3} V. Danilov, {\em De Rham complex on toroidal variety},
in {\sl Algebraic Geometry} (I. Dolgachev, W. Fulton eds.), Lecture
Notes in Math. {\bf 1479}, Springer-Verlag, Berlin-Heidelberg-New York,
1991, 26--38.

\bibitem{dan2} V. Danilov, {\em Newton polyhedra and vanishing
cohomology}, Funct. Anal. and App. {\bf 13} (1979), 103--115.

\bibitem{dan1} V. Danilov, {\em The geometry of toric varieties}.
Russian Math. Surveys {\bf 33} (1978), 97--154.

\bibitem{dan.hov} V. Danilov, A. Khovanski\^i, {\em Newton polyhedra
and an algorithm for computing  Hodge-Deligne numbers},
Math. USSR Izv. {\bf 29} (1987), 279--298.

\bibitem{deligne} P. Deligne, {\em Th\'eorie de Hodge}, II, III,
Publ. Math. I.H.E.S. {\bf 40} (1971), 5--58; {\bf 44} (1975), 5--77.

\bibitem{dem} M. Demazure, {\em Sous-groupes alg\'{e}briques de
rang maximum du grupe de Cremona}, Ann. Sci. \'{E}cole Norm.
Sup. {\bf 3} (1970), 507--588.

\bibitem{dolgach} I. Dolgachev, {\em Weighted projective varieties},
in {\sl Group Actions and Vector Fields} (J.B. Carrell ed.), Lecture Notes
in Math. {\bf 956}, Springer-Verlag, Berlin-Heidelberg-New York, 1982,
34--71.

\bibitem{griff1} P. Griffiths, {\em On the periods of certain rational
integrals} I, II, Ann. of Math. {\bf 90} (1969), 460--495, 498--541.

\bibitem{groth} A. Grothendieck, {\em On the de Rham cohomology of algebraic
varieties}, Publ. Math. I.H.E.S. {\bf 29} (1966), 95--103.

\bibitem{grunbaum} B. Gr\"unbaum, {\sl Convex Polytopes}, John Wiley \&
Sons, New York, 1967.

\bibitem{hart} R. Hartshorne, {\sl Algebraic Geometry}, Springer-Verlag,
Berlin-Heidelberg-New York, 1977.

\bibitem{ham} J. Humphreys. {\sl Linear algebraic groups}, Springer-Verlag,
Berlin-Heidelberg-New York, 1975.

\bibitem{kempf} G. Kempf, F.Knudsen, D. Mamford, B. Saint-Donat, {\em
Toroidal embeddings}, I, Lecture Notes in Math. {\bf 339},
Springer-Verlag, Berlin-Heidelberg-New York, 1973.

\bibitem{khov} A. Khovanski\^i, {\em Newton polyhedra and toroidal
varieties}, Funct. Anal. and App. {\bf 11} (1977), 289--296.

\bibitem{man.zfas} Yu. Manin, M. Tsfasman, {\em Rational varieties:
algebra, geometry, and arithmetics}, Russian Math. Surveys {\bf 41}
(1986), 51--116.

\bibitem{mam} D. Mumford, J. Fogarty, {\sl Geometric invariant theory}
($2^{\rm nd}$ ed.), Springer-Verlag, Berlin-Heidelberg-New York,
1982.

\bibitem{oda} T. Oda, {\em Convex Bodies and Algebraic Geometry},
Springer-Verlag, Berlin-Heidelberg-New York, 1988.

\bibitem{pet.steen} C. Peters, J. Steenbrink, {\em Infinitesimal variations of
Hodge structure and the generic Torelli problem for projective
hypersurfaces (after Carlson, Donagi, Green, Griffiths, Harris)\/}, in
{\sl Classification of Algebraic and Analytic Manifolds\/} (K.~Ueno
ed.), Progress in Math.~{\bf 39}, Birkh\"auser, Boston-Basel-Berlin,
1983, 399--463.

\bibitem{prill} D. Prill, {\em Local classification of quotients of
complex manifolds by discontinuous groups}, Duke Math. J. {\bf 34}
(1967), 375--386.

\bibitem{reisner} G. Reisner, {\em Cohen-Macaulay quotients of
polynomial rings}, Adv. in Math. {\bf 21} (1976), 30--49.

\bibitem{scheja} G. Scheja, U. Storch, {\em \"Uber Spurfunktionen bei
v\"ollst\"andigen Durchschnitten\/}, J. Reine u. Angew. Math. {\bf
278/279} (1975), 174--189.

\bibitem{stanley} R. Stanley, {\sl Combinatorics and Commutative
Algebra}, Progress in Math. {\bf 41}, Birk\-h\"auser,
Boston-Basel-Berlin, 1983.

\bibitem{steen} J. Steenbrink, {\em Intersection form for quasi-homogeneous
singularities}, Comp. Math. {\bf 34} (1977), 211--223.


\end{thebibliography}
\end{document}